\documentclass{amsart}[12pt]

\usepackage{geometry} 

\def\Ree{\mathbb{R}}

\def\Nee{\mathbb{N}}

\usepackage{hyperref} 
\usepackage{stmaryrd}

\begin{document}

\title{Compactness Properties for Geometric Fourth Order Elliptic Equations with Application to the $Q$-curvature flow }

\author{Ali Fardoun}
\address{Laboratoire de Math\'ematiques, UMR 6205 CNRS 
Universit\'e de Bretagne Occidentale 
6 Avenue Le Gorgeu, 29238 Brest Cedex 3   
France}
\email{Ali.Fardoun@univ-brest.fr}
\author{Rachid Regbaoui }
\address{Laboratoire de Math\'ematiques, UMR 6205 CNRS 
Universit\'e de Bretagne Occidentale 
6 Avenue Le Gorgeu, 29238 Brest Cedex 3   
France}
\email{Rachid.Regbaoui@univ-brest.fr}


\begin{abstract}

We prove the compactness of solutions to general  fourth order elliptic equations which are  $L^1$-perturbations of the $Q$-curvature equation  on   compact Riemannian  $4$-maniflods. Consequently, we  prove the global existence and convergence of the $Q$-curvature flow on a generic class of Riemannian $4$-manifolds. As a by product,  we give a positive answer to an open  question by A. Malchiodi \cite{aM1}  on the existence of bounded Palais-Smale sequences for the $Q$-curvature problem when the Paneitz operator is positive with  trivial kernel. 

\end{abstract}

\maketitle

\section{Introduction and statement of the results }

\medskip

On a Riemannian $4$-manifold  $(M, g_0)$,  the $Q$-curvature  of the metric $g_0$  is defined by 

\medskip

$$Q_0 = -\frac{1}{6}\left( \Delta_0S_0 + S_0{}^2 - 3| {\rm Ric}_0|^2\right) , \eqno (1.1)
$$

\bigskip

\noindent where $\Delta_0$, $S_0$ and ${\rm Ric}_0$ denote respectively, the
Laplace-Beltrami operator,
the scalar curvature and the Ricci tensor associated to the metric $g_0$.
A conformal change of the metric $g_0$ produces a metric $g = e^{2u}g_0$
having $Q$-curvature

\medskip

$$
Q_g = e^{-4u}(P_0u + Q_0),  \eqno (1.2)
$$

\bigskip

\noindent where $P_0$ is the Paneitz operator defined by

\medskip

$$
P_0u = \Delta_0^2 u + {\rm div}_0\left[\left( {2\over 3}S_0g_0 - 2{\rm
Ric}_0\right)du\right],   \eqno (1.3)
$$

\medskip

\noindent where  ${\rm div}_0$ denotes the divergence operator with respect to $g_0$. The Paneitz operator is conformally invariant in the sense that if   $g = e^{2u}g_0$,  then  the Paneitz operator with respect to $g$ is given by $ P_g = e^{-4u} P_0$.   Through this paper,   we will always assume that $P_0$ has  trivial kernel, that is, its kernel consists only of constant functions. Some times, we will  also  need to assume that $P_0$ is positive, which means that for all $u \in C^{\infty}(M)$ we have 
$$\int_M P_0u\cdot u \ dV_0 \ge 0,$$

\medskip

\noindent where $dV_0$ is the volume element with respect to $g_0$.  We note here that both hypothesis are conformally invariant.  From the following  Chern-Gauss-Bonnet formula :

$$
\int_MQ_g\, dV_g + \frac{1}{4}\int_M |W_g|^2 \, dV_g = 8\pi^2 \chi (M)\,,
$$

\medskip

\noindent where $\chi (M)$ denotes the Euler-characteristic of $M$,  $W_g$ is the  Weyl  tensor of $g$  and $dV_g$ is the volume element with respect to $g$, we see that the total $Q$-curvature 

\medskip

$$
k_0 := \int_MQ_0\, dV_0 =  \int_MQ_g\, dV_g  \eqno (1.4)
$$

\medskip

\noindent is also conformally invariant since the Weyl tensor is pointwise conformally invariant.  We note here that formula (1.4) is also a direct  consequence of (1.2). 

\medskip

The $Q$-curvature and the Paneitz operator have
received much attention in recent years because of
their role in four-dimensional conformal geometry and in
mathematical physics. Similarly to the uniformization problem of surfaces, one of the interesting problems in the geometry of 4-manifolds is to ask if there exists a  metric $g$ on $M$  conformal to $g_0$ and having constant  $Q$-curvature ?  By using the relation (1.2), the problem is equivalent to find a function $u \in C^{\infty}(M)$ satisfying the following partial  differential equation 

\medskip

$$P_0 u + Q_0 = k_0 e^{4u}. \eqno (1.5) $$

\medskip 

Equation (1.5) has a variational structure, and its solutions are critical points of the following functional on the Sobolev space $H^2(M)$ :

\medskip

$$E(u) := {1\over 2} \int_{M} P_0 u\cdot u \ dV_0 +  \int_M Q_0 u \ dV_0 - {k_0 \over 4} \log\left(\int_M e^{4u} dV_0\right)  ,   \eqno (1.6) $$

\medskip
 
 \noindent that we will call the $Q$-curvature  functional through this paper.  The  main difficulty in the study of this functional is that in general it is not coercive, and it can be  unbounded from below and from above. This is due   to the large values that  the total $Q$-curvature $k_0$ may have, and to the possibility of negative eigenvalues of the operator $P_0$. 
 
 \medskip 
 
 S.A. Chang  and  P. Yang \cite{aC} first studied Equation (1.5)   by minimizing the  functional $E$. They constructed conformal metrics of constant $Q$-curvature when  the Paneitz operator $P_0$ is positive with  trivial kernel and the total $Q$-curvature satisfies $k_0 < 16\pi^2$ which  is  the total $Q$-curvature of the Euclidean sphere $\mathbb{S}^4$. The key point in their proof is that   by using Adams inequality  (see section 2),  if we suppose that $P_0$ is positive with trivial kernel and $k_0 < 16 \pi^2$, then the functional $E$ is bounded from below, coercive and its critical points can be found as global minima. Later,  Z.Djadli and A.Malchiodi \cite{zD}  solved  equation (1.5) when $P_0$ is not necessarily positive under the condition that the kernel of $P_0$ is trivial and  $k_0 \not = 16 k \pi^2, \  k \in \mathbb{N}^*$.  They constructed a critical point of $E$ by a mini-max scheme  based on a result by  A. Malchiodi \cite{aM1}, and independently   by  O. Druet and F. Robert \cite{oD},  on the compactness of solutions to fourth 
 order elliptic equations.  For more details on the $Q$-curvature problem we refer the reader to \cite{pB1}, \cite{pB2}, \cite{sB2}, \cite{mG} and the  references therein.

 \medskip

 In \cite{aM1} A. Malchiodi  proved, by assuming   $k_0 \not\in 16  \pi^2  \mathbb{N}^*,$   the compactness of  any sequence $(u_n)_n$ satisfying a $C^0$-perturbation of equation (1.5) of the form
 
 \medskip

$$P_0u_n + Q_n = k_ne^{4u_n}  \eqno (1.7)$$

\noindent where $\displaystyle k_n  = \int_M Q_n dV_0$, by assuming 

 $$Q_n \underset{n \to + \infty}{\longrightarrow} Q_0 \  \ \hbox{in} \    C^{0}(M) .$$
 
 \medskip
 
This result does not apply to Palais-Smale sequences for the functional $E$ since for such sequences one needs $H^{-2}$-perturbations of equation (1.5).  An open question was kept in \cite{aM1}:  do there exist    bounded Palais-Smale sequences for  the functional $E$ ? One of the main result of the present paper is to give a positive answer to this question  when the Paneitz operator is positive with trivial kernel. We prove  first a compactness result  of solutions to $L^1$-perturbations of equation (1.5), that is,   we need only to assume  
 
 $$Q_n \underset{n \to + \infty}{\longrightarrow} Q_0 \  \ \hbox{in} \    L^1(M) $$
 
 \noindent in equation (1.7) above. Then we apply this result to study  the solution of  the heat flow equation associated with (1.5) since such an equation is a  $L^1$-perturbation of (1.5),  giving  thus a convergent Palais-Smale  sequence for the  functional $E$.  Our first main result reads :
  
  \medskip
  
 \newtheorem{theo}{Theorem}[section]

    \begin{theo}Let $(M,g_0)$ be a compact Riemannian $4$-manifold  whose  Paneitz  operator  has  trivial kernel.  Let $f \in C^0(M)$, and  let $(u_n, f_n)_n$ be a sequence in $H^2(M)\times L^1(M)$ satisfying

$$P_0u_n + f_n = k_ne^{4u_n}    \eqno (1.8)$$

$$  \int_M e^{4u_n} dV_0 = 1 \eqno (1.9) $$
and
 $$f_n \underset{n \to + \infty}{\longrightarrow} f \  \ \hbox{in} \    L^1(M)  \  ,  \eqno (1.10)$$ 
 
\noindent  where  \  $\displaystyle k_n = \int_Mf_n dV_0$.  Then  one   of the following alternatives holds :
 
 \bigskip

   \noindent $1)$ either the  sequence $(e^{|u_n|})_n$ is bounded in $L^{p}(M)$ for all $p \in [1, + \infty)$, 
   
   \bigskip
   
   \noindent $2)$ or for a subsequence, that we still denote  by $(u_n)_n$ for simplicity, there exist a finite number of points $a_1, ... , a_m \in M$ and integers $l_1, ... , l_m \in \Nee^*$ such that 
   
   \medskip
   
   $$\sum_{j=1}^ml_j  =   {1 \over 16\pi^2}  \int_M fdV_0 \eqno(1.11)$$
   and 
   $$e^{4u_n} \underset{n \to + \infty}{\longrightarrow} {16\pi^2\over \int_{M}fdV_0}  \sum_{j=1}^ml_j \delta_{a_j}, \eqno (1.12)$$
in the sense of measures, where $\delta_a$ stands for the Dirac mass at the point $a \in M$.

\end{theo}

   \bigskip

A particular case of Theorem 1.1 is the following result when  $ \displaystyle \int_Mf dV_0  \not\in 16\pi^2\Nee^* $. 
 
 \medskip
 
  \newtheorem{cor}{Corollary}[section]
 
     \begin{cor}Let $(u_n, f_n)_n$ be a sequence in $H^2(M)\times L^1(M)$ as in Theorem 1.1. If  we assume in addition
 $$\int_Mf dV_0  \not\in 16\pi^2\Nee^*,     \eqno (1.13)$$
\medskip
 \noindent then   the  sequence $(e^{|u_n|})_n$ is bounded in $L^{p}(M)$ for all $p \in [1, + \infty)$, that is,
 
 $$\int_Me^{p|u_n|} dV_0 \le C_p , $$
 
 \medskip
 
 \noindent where $C_p$ is a positive constant depending on $p$ but not on $n$.

\end{cor}
 
\bigskip

\medskip

It is clear that Corollary 1.1  is a particular case of Theorem 1.1 since if  condition (1.13) is satisfied, then the second alternative in Theorem 1.1 does  not occur.  

\medskip

As a direct consequence of Corollary 1.1,  we recover the result of A. Malchiodi \cite{aM1} and O. Druet-F. Robert \cite{oD} stated above. More precisely, we have  the following corollary  on the compactness of solutions to the $Q$-curvature equation (1.5) : 

\medskip

     \begin{cor}Let $(M,g_0)$ be a compact Riemannian $4$-manifold    whose Paneitz  operator has a  trivial kernel, and assume that the total $Q$-curvature $k_0$ satisfies  $k_0  \not\in 16\pi^2\Nee^*$.   Then for any \  $ k \in \Nee$,    there  exits a constant $C_k$ depending only on $k$ and  $(M, g_0)$ such that  for any solution $u \in H^2(M)$ of   the $Q$-curvature equation 
     $$ P_0 u + Q_0 = k_0e^{4u}$$
     with the normalization  \ $\displaystyle  \int_M e^{4u} dV_0 = 1$,  we have 
     $$\|u\|_{C^k(M)} \le C_k .$$
     
\end{cor}

\bigskip

We note here that when $k_0 \not=0$, the normalization condition  $\displaystyle  \int_M e^{4u_n} dV_0 = 1$ is automatically satisfied from  the $Q$-curvature equation.

\medskip

\bigskip

 \newtheorem{rem}{Remark}[section]

\begin{rem} 1) The conclusions of Theorem 1.1 and Corollary 1.1  remain  valid if we consider equations of the form 
$$ P_0u_n + f_n =  h_n e^{4u_n} $$
where we assume $  h_n \in C^0(M)$ such that $C^{-1} \le h_n \le C$ for some positive constant $C$ independent of $n$.  This can be checked by following the same proofs with some necessary slight modifications.  
\medskip

\noindent 2) One can easily check that in Corollary 1.1 a subsequence of $(u_n)_n$ converges strongly in $W^{2, p}(M)$ for all $p \in [1, 2)$ to a function $u_{\infty}$ satisfying the following equation 
$$P_0 u_{\infty} + f = k_0e^{4u_{\infty}} ,$$
where $\displaystyle  k_0 = \int_M f dV_0$. 

\noindent 3) As it can be seen from our proofs, Theorem 1.1 and Corollary 1.1 remain valid if we consider $L^1$-perturbations of the mean field equation on compact Riemannian surfaces. Indeed, the same arguments work if we replace the  Paneitz operator by the Laplacian. Some related results  concerning the mean field equation are proved by J-B. Cast\'eras in his thesis \cite{jC}. 
 
 \end{rem}

\bigskip

The main difficulty in the study of equations like (1.8) in Theorem 1.1  is the appearance of the so-called bubbling phenomena due to the concentration of the volume of the conformal metric $g_n = e^{2u_n}g_0$.  We prove that if such  phenomena  occur then there must  be some volume quantization. An important tool in the proof of such a result is an integral Harnack type inequality that we will prove in section 3.

\bigskip

As stated above, the second main result of this paper  concerns  the evolution problem associated with equation (1.5). More precisely,  we  will   consider  the evolution of a metric $g$ on $M$  under the 
flow:

\medskip

$$
\begin{cases} \partial_t g  = - \left( Q_g -  \overline{Q}_g   \right) g \cr \cr
 g(0)= e^{2u_0} g_0,  \    u_0 \in C^{\infty} (M) \end{cases}  \eqno (1.14)
$$ 
\medskip

\noindent where
 $$  \overline{Q}_g  = {1 \over \hbox{Vol}_g(M)}\int_MQ_g dV_g =  {k_0 \over \hbox{Vol}_g(M)}$$
 
 \medskip
 
 \noindent  is  the average of $Q_g$. 
 
 \medskip

 Since equation (1.14) preserves the conformal structure of $M$, then \ $ g(t) = e^{2u(t)} g_0$, where $u(t) \in C^{\infty}(M)$ with initial condition $u(0) = u_0  \in C^{\infty}(M)$. For simplicity,  we have used the notation  $u(t):= u(. , t),  \ t \in I$,  for any function defined on $M\times I$,  where $I$ is a subset of $\Ree$. Thus the flow (1.14) takes the form

\medskip

 $$ \displaystyle 
\begin{cases} \partial_t u  = -{1 \over 2} e^{-4u} \left( P_0u + Q_0 \right) + {1 \over 2}\frac{ k_0 }{\int_M  e^{4u} \ dV_0 }, \cr  \cr u(0)= u_0.   \end{cases}
\eqno (1.15)
$$

\bigskip

\medskip

It is clear that  the first equation in  (1.15) is parabolic since $P_0$ is an elliptic operator. Then by classical methods it admits a solution $ u \in C^{\infty} (M\times [0,T) )$ where $T \leq + \infty$ denotes the maximal time of existence. By integrating  the first equation in (1.15) over $M$ with respect to the volume element of $g(t)$, we see that the volume of $M$ with respect to $g(t)$ remains constant, that is, 
$$\int_M e^{4u(t)} dV_0 =  \int_M e^{4u_0} dV_0 ,  \  \   \forall t \in [0, T). \eqno (1.16)$$

\medskip

If we multiply the first equation in  (1.15) by $\partial_tu$ and integrating with repect to $g(t)$, we see that the functional $E$ is decreasing along the flow :
$${d\over dt}E(u(t)) = - 2 \int_M e^{4u(t)} |\partial_tu(t)|^2 dV_0,   \  \   \forall t \in [0, T).  \eqno (1.17) $$

\bigskip

As far as we know, the evolution problem (1.15) has been studied  only in the case where the total $Q$-curvature $k_0$ satisfies $k_0 \le 16\pi^2$ and $P_0$  is positive with trivial kernel. Indeed,  S.Brendle \cite{sB1} was the first who studied the $Q$-curvature flow by  considering a more general flow (with prescribed  $Q$-curvature function) on 
 Riemannian  manifolds of even dimension. In dimension four,  his result corresponds to assume  that  $P_0$ is positive with  trivial kernel and $k_0<  16 \pi^2$, and then  he proved that  (1.15)  has
a solution which is defined for all time $(T= + \infty)$ and converges to a smooth function $u_{\infty}$ such that the metric $g_{\infty} = e^{2u_{\infty}}$ has  constant $Q$-curvature.  When $k_0 =  16 \pi^2$, he proved in \cite{sB3} the  global existence and the convergence of $Q$-curvature flow on the sphere $\mathbb{S}^4$. We also mention here the work of A. Malchiodi and M. Struwe \cite{aM2} where they consider the $Q$-curvature flow on $\mathbb{S}^4$  with a prescribed $Q$-curvature function $f$. They proved the global existence of the flow and studied its asymptotic behaviour  under some assumptions on the critical points of $f$.

\medskip

Our second main result in this paper  is to study  the flow  (1.15)  on  Riemannian 4-manifolds   with total curvature $k_0$ satisfying  $k_0 \not\in 16\pi^2\Nee^*$.  In particular,  we are able to allow  $k_0$ to  take values beyond the critical threshold $16\pi^2$.  
Our result is as follows : 

\medskip

\begin{theo}Let $(M,g_0)$ be a compact Riemannian $4$-manifold whose  Paneitz  operator $P_0$   is positive with  trivial kernel.  For $u_0 \in  C^{\infty}(M)$, let $u \in C^{\infty}(M\times [0, T))$   the solution of problem (1.15) defined on a maximal interval $[0, T)$.  If

$$\inf_{t \in [0, T)} E(u(t)) > - \infty \   ,   \eqno (1.18)$$

\noindent where $E$ is defined  by (1.6),    then  $T = + \infty$, that is, $u(t)$ is globally defined on $[0, + \infty)$. Moreover,  if in addition the total $Q$-curvature $k_0$ satisfies $  k_0 \not\in 16\pi^2\Nee^*$,  
then  $u(t)$   converges in $C^{\infty}(M)$ as  $t\to + \infty$,  to a function $u_{\infty} \in C^{\infty}(M)$ satisfying  the $Q$-curvature equation 

$$ P_0 u_{\infty} + Q_0 ={ k_0 \over \int_Me^{4u_{\infty}} dV_0}  e^{4u_{\infty}} . $$

\end{theo}

\bigskip

 It is natural to ask if  there exist initial data  $u_0 \in C^{\infty}(M)$ for which the  solution $u(t)$  of problem (1.15)  with $u(0) = u_0$, satisfies condition (1.18)  in Theorem 1.2 ?  As we will see below,  we can always find initial data  for which (1.18)  is satisfied,  and others   for which it is not the case.  We note  here that the condition (1.18) is  automatically  satisfied if $ k_0 \le 16 \pi^2$  and $P_0$  positive with trivial kernel since in this case  the functional $E$ is bounded from below by using Adams inequality (see section 2).   First we have : 

\medskip

\begin{theo}Let $(M,g_0)$ be a compact Riemannian $4$-manifold whose  Paneitz  operator $P_0$ is positive with  trivial kernel, and  suppose that the total $Q$-curvature $k_0$ satisfies $  k_0 \not\in 16\pi^2 \mathbb{N}^*$.   Then there exists at least one function $u_0 \in  C^{\infty}(M)$ such that the solution $u( t)$ of problem (1.15)  with $u(0) = u_0$, satisfies condition (1.18) in  Theorem 1.2,  that is, 

$$\inf_{t \in [0, T)} E(u(t)) > - \infty \   .$$

\medskip

\noindent Thus, according to Theorem 1.2,   $u(t)$  is globally defined on $[0, + \infty)$ and   converges in $C^{\infty}(M)$ as  $t\to + \infty$,  to a function $u_{\infty} \in C^{\infty}(M)$ satisfying  the $Q$-curvature equation 
$$ P_0 u_{\infty} + Q_0 ={ k_0 \over \int_Me^{4u_{\infty}} dV_0}  e^{4u_{\infty}} . $$

\end{theo}

\bigskip

It follows from  Theorem 1.3, by taking any  real sequence  $(t_n)_n$ such that \  $ t_n \underset{n \to + \infty}{\longrightarrow}  + \infty$,  that $(u(t_n))_n$ is a convergent Palais-Smale sequence for the $Q$-curvature functional $E$. This gives  a positive answer to the open  question by  A. Malchiodi \cite{aM1} stated above.  In particular,  Theorem 1.3  is a direct method to solve equation (1.5).  

\medskip

\medskip

The following theorem gives a class of functions $u_0 \in C^{\infty}(M)$ for which condition (1.18) in Theorem 1.2 is not satisfied. More precisely, it gives a class of initial data for which the  corresponding flow  $u( t)$ blows up in finite or infinite time :

\medskip

\begin{theo}Let $(M,g_0)$ be a compact Riemannian $4$-manifold whose  Paneitz  operator $P_0$ is positive with  trivial kernel, and  suppose  that the total $Q$-curvature $k_0$ satisfies $  k_0 \not\in 16\pi^2 \mathbb{N}^*$.   Then there exists  a constant $\lambda \in \Ree$ such that for any $u_0 \in C^{\infty}(M)$ satisfying  $E(u_0) \le \lambda$, the corresponding solution $u( t)$ of  (1.15) satisfies  $ \displaystyle  \lim_{t\to T} E(u( t)) = - \infty .$
 \end{theo}

\bigskip

One of the principal  difficulties in the study of the $Q$-curvature flow (1.15) comes from the absence of a maximum principle for elliptic operator of high order (greater than $4$). Our analysis is based on the proof of some delicate  integral estimates  concerning the parabolic equation (1.15), combined with the compactness of the solutions   of the corresponding elliptic equation  proved in Theorem 1.1.

\medskip

\section{Preliminaries and blow-up analysis}

\medskip

We introduce in this section some basic  tools on elliptic operators on Riemannian manifolds, and we recall some known results  on the blow-up of solutions to general $Q$-curvature  type equations.  

\medskip

Let $(M, g_0)$ a smooth compact 4-Riemannian manifold without boundary. For the simplicity of notations,  the Riemannian distance  between two points $x, y \in M$ is denoted as in the Euclidean space, by $|x-y|$. If $x \in M$ and $r>0$, we denote by $B_r(x)$ the geodesic ball in $M$ of center $x$ and radius $r$.  If  $ x \in \Ree^n$,  we denote in the same manner by  $B_r(x)$ the Euclidean ball in $\Ree^n$ of center $x$ and radius $r$.  The volume element of $g_0$ is denoted by $dV_0$, and the volume of  any measurable set $A\subset M$ is denoted by $|A|$.  The Green  function $G$  associated to  the Paneitz  $P_0$ is a symmetric function  $G \in C^{\infty}(M \times M \setminus D )$, 
where $D= \{ (x, x) \ : \ x \in M \}$ is the diagonal of $M$, giving the inversion formula  for Paneitz operator. That is, if $F\in L^1(M)$ with $\overline{F} = 0$, then $u$ is a solution of 

$$P_0 u = F \eqno (2.1) $$
if and only if 
$$u(x) = {\bar u} + \int_MG(x,y) F(y) dV_0(y) , \eqno (2.2)$$
where we denote by ${\bar h} := {1 \over |M|}\int_M h dV_0$ the average of any function $h \in L^1(M)$. We have the following asymptotics for $G$ 
$$G(x, y) = - {1 \over 8 \pi^2} \log|x-y| +  R(x,y), \eqno (2.3)$$
where $R \in C^{0}(M\times M) $.  

\medskip

The following proposition concerning solutions of equation (2.1) is proved in A. Malchiodi \cite{aM1} by using the asymptotics of the Green function (see Lemma 2.3 in \cite{aM1}.)

\medskip

\newtheorem{prop}{Proposition}[section]

\begin{prop} Let $(u_n, F_n)   \in H^2(M)\times L^1(M)$ satisfying 
$$P_0u_n = F_n$$
with $\|F_n \|_{L^1(M)} \le K$ for some constant $K$ independent of $n$. Then for any  $x \in M$, for any  $r >0$ (small enough),  for  any $j= 1, 2, 3 $,  and $p \in [1, 4/j)$, we have 
$$\int_{B_r(x)}|\nabla^ju_n|^p dV_0 \le  Cr^{4-jp} , $$ 
where $C$ is a positive constant depending on $K, M, p$ but not on $n$. 
\end{prop}

\bigskip

We need also the following proposition proved in \cite{aM1} :

\medskip

\begin{prop}  Let $(u_n, F_n)  \in H^2(M)\times L^1(M)$ satisfying 
$$P_0u_n = F_n$$
with $\|F_n \|_{L^1(M)} \le K$ for some constant $K$ independent of $n$. Then  :

\medskip

\medskip

\noindent 1) either 
$$\int_M e^{q (u_n - {\bar u_n})} dV_0 \le C$$
for some  $q>4$ and some $C >0$ (independent of $n$), 

\medskip

\medskip

\noindent 2) or there exists  a point $x \in  M$ such that for any $r >0$, we have 
$$\liminf_{n\to +\infty} \int_{B_r(x)}|F_n| dV_0 \ge 8 \pi^2 .$$
\end{prop} 

\bigskip

\medskip

\begin{rem} Proposition 2.2  remains valid if one replace the metric $g_0$ on $M$ by a family of metric $(g_n)_n$ depending on $n$ which is uniformly bounded in $C^k(M)$ for any $k \in \mathbb{N}$. The same result holds also if we replace $M$ by any bounded open ball of $\mathbb{R}^4$ and  assuming all the functions with compact support in this ball. Indeed,  a bounded open ball  of $\mathbb{R}^4$ can  always be embedded in a  torus for example. 
\end{rem}

\bigskip

Now we  shall give some basic properties of solutions to equation (2.1) when $F $ is as in Theorem 1.1.  That is, we consider  a sequence $(u_n, f_n)_n$ in    $H^2(M)\times L^1(M)$ satisfying 

\medskip

$$P_0u_n + f_n = k_ne^{4u_n} ,   \eqno (2.4)$$

\medskip

\noindent such that   
 $$ \int_M e^{4u_n} dV_0 = 1,  \eqno (2.5) $$
 and 
 $$f_n \underset{n \to + \infty}{\longrightarrow} f \  \ \hbox{in} \    L^1(M)  \  ,  \eqno (2.6 )$$ 
 
 \noindent with $f \in C^0(M)$, and where $\displaystyle k_n = \int_Mf_n dV_0$.  
 
\bigskip

First we state the following proposition which can easily be deduced from Proposition 2.2 above by setting $F_n= k_ne^{4u_n} - f_n$. 

\medskip

\begin{prop} Let $(u_n, f_n)_n  \in   H^2(M)\times L^1(M)$ satisfying (2.4)-(2.6). Then  :

\medskip

\medskip

\noindent 1) either  for any $p \ge 1$ we have for some   some constant $C_p >0$ independent of $n$, 
$$\int_M e^{p|u_n |} dV_0 \le C_p$$

\medskip

\noindent 2) or there exists  a point $x \in  M$ such that for any $r >0$, we have 
$$k_0 \liminf_{n\to +\infty} \int_{B_r(x)}e^{4u_n}dV_0 \ge  8 \pi^2 + o_r(1) ,$$
where $\displaystyle k_0 = \int_M fdV_0$,  and where $o_r(1) \to 0 $ as $r \to 0$.
\end{prop} 

\bigskip

\medskip

  If  $(u_n)_n$ and $(f_n)_n$ are as above, $(x_n)_n$ is a sequence of points in $M$ and $(r_n)_n$ a sequence of positive numbers such that $r_n \to 0$, we set 
 $$\widehat {u}_n(z) = u_n(\exp_{x_n}(r_n z) ) + \log r_n , \  \   z \in B_{\delta\over r_n}(0) , \eqno (2.7) $$
 where $ B_{\delta\over r_n}(0) \subset \Ree^4$ is the Euclidean ball of center $0$ and radius ${\delta\over r_n}$, and where $\delta$ is the injectivity radius of $M$.  We note here that  the ball $ B_{\delta\over r_n}(0)$ approaches $\Ree^4$ when $n \to + \infty$.  As we will see later, it is useful to introduce  the following quantities. Let  $T_n :   B_{\delta\over r_n}(0) \to M$ defined by $T_n(z) = \exp_{x_n}(r_n z)$,  and define a  metric $g_n$ on 
 $B_{\delta\over r_n}(0)$  by 
 $$g_n = r_n^{-2}T_n^{*}g.  \eqno (2.8)$$

\medskip

\noindent It is not difficult to see that $g_n  \underset{n \to + \infty}{\longrightarrow} g_{\mathbb{R}^4}$  in $C^{k}(B_R(0))$ \  for all $k \in \mathbb{N}$ and all $R>0$, where $g_{\mathbb{R}^4}$  is the standard Euclidean metric of $\mathbb{R}^4$.  An easy computation shows that $\widehat{u}_n$ satisfies the following PDE in $\mathbb{R}^4$ 

$$P_{g_n}\widehat{u}_n + r_n^4  \widehat{f}_n = k_ne^{4\widehat{u}_n} , \eqno (2.9) $$

\medskip

\noindent where $P_{g_n}$ is the Paneitz operator of the metric $g_n$ in $\mathbb{R}^4$  and  

$$\widehat{f}_n(z) = f_n(\exp_{x_n}(r_n z) ), \   z \in B_{\delta\over r_n}(0).  \eqno (2.10) $$ 

\medskip

We shall also use the following function on $\Ree^4$, known as a standard bubble, 

 $$\xi_{z_0}(z) = \log\left( {2\lambda  \over 1 +  \lambda^2  |z-z_0|^2 }\right) - {1\over 4}\log\left({k_0 \over 6}\right),  \ z \in \Ree^4, \eqno (2.11) $$
 
 \noindent where $z_0$ is  a fixed point in $\Ree^4$,  \ $\lambda > 0$ is  a positive constant, and $\displaystyle k_0 =    \int_M f dV_0$ that we assume satisfying $k_0 > 0$. 
 
 \medskip

 The following is a slightly different definition of blow-up  with respect to  that given  in \cite{aM1} since we are considering more general equations. 
 
 \medskip
 
 \newtheorem{Def}{Definition}[section]
 
 \begin{Def} Let $(u_n, f_n)_n  \in   H^2(M)\times L^1(M)$ satisfying (2.4)-(2.5). Let $(x_n)_n$ a sequence of points in $M$ and $(r_n)_n$ a sequence of positive numbers such that $r_n \to 0$. We say that the sequence $(x_n, r_n)_n$ is a blow-up for $(u_n)_n$  if  for some $z_0 \in \Ree^4$, we have  for any $R >0$, 
 $$\widehat{u}_n  \underset{n \to + \infty}{\longrightarrow}  \xi_{z_0}  \hspace{2mm}   \hbox{in}  \hspace{2mm}    {\mathcal D}^{\prime}(B_R(0)) \hspace{3mm}   \hbox{and}     \   \    e^{4\widehat {u}_n} 
 \underset{n \to + \infty}{\longrightarrow}    e^{4\xi_{z_0}}   \    \   \hbox{in} \  \  L^1(B_R(0)),  \eqno (2.12)$$
where $\widehat {u}_n$ and $\xi_{z_0}$  are defined by (2.7) and (2.11), and where  ${\mathcal D}^{\prime}(B_R(0))$ denotes the space of distributions on $B_R(0)$. 
 \end{Def}
 
\medskip

By using the above definition and the fact that $\displaystyle \int_{\Ree^4} e^{4\xi_{z_0}} dz = {16\pi^2 \over k_0}$, one can easily prove the following : 

\medskip

\begin{prop} Let $(u_n, f_n)_n  \in   H^2(M)\times L^1(M)$ satisfying (2.4)-(2.5), and let  $(x_n, r_n)_n$  a blow-up for $(u_n)_n$. Then for any positive  sequence $(\beta_n)_n$ such that $\beta_n \to +\infty$, there exists a positive sequence $(b_n)_n$ with  $b_n \le \beta_n $ and $b_n \to +\infty$ such that 
$$\lim_{n\to + \infty} \int_{B_{b_nr_n}(x_n)} e^{ 4u_n} dV_0 = {16\pi^2 \over k_0} .$$
\end{prop}

\medskip

Now we will prove a proposition giving the existence of blow-ups for solutions to equation (2.4).

\medskip

\begin{prop} Let $(u_n, f_n)_n  \in   H^2(M)\times L^1(M)$ satisfying (2.4)-(2.5). Suppose that there exit a sequence $(x_n)_n$ in $M$,  sequences $(r_n)_n, (\widehat{r}_n)_n$  of positive numbers and a constant $\rho \in (0, {\pi^2\over k_0} ]$  (independent of $n$ ) such that 
$$\widehat{r}_n \underset{n \to + \infty}{\longrightarrow}  0 \ \ , \  \  {r_n \over \widehat{r}_n}  \underset{n \to + \infty}{\longrightarrow}  0  \eqno (2.13) $$
and  for  $n$ large enough, 
$$  \int_{B_{r_n}(x_n)}e^{4u_n} dV_0 = \rho \ \ , \  \  \int_{B_{r_n}(y)}e^{4u_n} dV_0 \le {\pi^2\over k_0}  \   \    \forall \ y \in B_{\widehat{r}_n}(x_n) . \eqno (2.14) $$

\medskip

\noindent Then by passing to a subsequence,  $(x_n, r_n)_n$ is a blow-up for   $(u_n)$. 
\end{prop}

\medskip

\begin{proof} The proof follows closely that of Proposition 3.4 in A. Malchiodi \cite{aM1}.  But since we are considering more general equations, we shall give the detailed  proof.  It suffices to prove that for any fixed 
$R \ge 1$,  a subsequence of  $(\widehat{u}_n)_n$  converges in $W^{2, p}(B_R(0))$   to $\xi_{z_0}$ for some $z_0 \in \mathbb{R}^4$ and some $p \ge 1$,   and  that $(e^{4\widehat{u}_n})_n$ converges  to $e^{4\xi_{z_0}}$  in  $L^1(B_R(0))$. 

\medskip

   Let  $R \ge 1$  be fixed  large enough.  In what follows,  $C$ denotes a positive constant depending on $M$ and $R$ but  independent of $n$, whose values may change from line to line. Set  

$$a_n  = {1 \over |B_{2R}(0)|_{g_n}} \int_{B_{2R}(0)} \widehat{u}_n dV_n,  $$

\medskip

\noindent where $dV_n$ is the volume element of the metric $g_n$ defined in (2.8) above, and $|B_{2R}(0)|_{g_n}$ is the volume of the ball $B_{2R}(0)$ with respect to the metric $g_n$. Let now $\varphi \in C^{\infty}_0(B_{2R}(0))$ such that $\varphi = 1 $ on $B_R(0)$, and define $v_n:= \varphi(\widehat{u}_n- a_n )$.  Then by using  Proposition 2.1 and the Poincar\'e inequality, we have  for any $p \in [1, 4/3)$ and $j = 1,2, 3$, 
$$\int_{B_{2R}(0)}|\nabla^j \widehat{u}_n|^p dV_n \le C \  \  \hbox{and} \  \  \int_{B_{2R}(0)} |\widehat{u}_n - a_n|^p dV_n \le C,$$
 which implies that  $(v_n)_n$   is bounded in $W^{3, p}(B_{2R}(0))$ for  any $p \in [1, 4/3)$, that is
$$ \|v_n\|_{W^{3,p}(B_{2R}(0))} \le C . \eqno (2.15) $$
(we recall here that the metric $g_n$ is uniformly bounded in $C^k(B_R(0))$ for any $k \in \mathbb{N}$ and $R> 0$.)

\medskip

\noindent Moreover we have 

$$P_{g_n}v_n = - r_n^4 \varphi \widehat{f}_n + k_n\varphi e^{4\widehat{u}_n}  +  h_n , \eqno (2.16) $$

\medskip

\noindent where $ h_n = L_n(\widehat{u}_n-a_n)$ for some third order linear operator $L_n$ with uniformly bounded smooth coefficients, and where  $\widehat{f_n}$ is defined by (2.10). So, $h_n$  is bounded in $L^p(B_{2R}(0))$ for any $p  \in [1, 4/3)$. If we set $F_n = - r_n^4 \varphi \widehat{f}_n + k_n\varphi e^{4\widehat{u}_n}  +  h_n$, then by using (2.6),  (2.13)-(2.14) and the fact that $(h_n)_n$ is bounded in $L^p(B_{2R}(0))$,  one can check that  there exists a constant $\delta_0 >0$ independent of $n$ such that for any $x \in B_{2R}(0)$  we have 
$$\int_{B_{\delta_0}(x)} |F_n| dV_n \le  2 \pi^2  + o_n(1),  \eqno (2.17) $$
where $o_n(1) \to 0$ as $n \to +\infty$. Thus by applying Proposition 2.2 (with Remark 2.1) to equation (2.16), we obtain that 
$$\int_{B_{2R}(0)}e^{q(v_n - \overline{v}_n)} dV_n \le C \eqno (2.18) $$
for some  constant  $q > 4$  independent of $n$, where $\overline{v}_n = \displaystyle {1 \over |B_{2R}(0)|_{g_n}} \int_{B_{2R}(0)} v_n dV_n$\  is the average of $v_n$ on $B_{2R}(0)$. Since $|\overline{v}_n| \le C$ by (2.15), then (2.18) becomes

$$\int_{B_{2R}(0)}e^{qv_n} dV_n \le C,   $$
which gives 
$$\int_{B_{R}(0)}e^{q(\widehat{u}_n-a_n)} dV_n \le C,   \eqno (2.19) $$
since $v_n = \widehat{u}_n-a_n$ on $B_R(0)$.

\medskip

 By using Jensen inequality and the fact 
 
$$\int_{B_{2R}(0)} e^{4\widehat{u}_n} dV_n  = \int_{B_{2r_nR}(x_n)} e^{4u_n} dV_0 \le \int_M e^{4u_n} dV_0 = 1,  \eqno (2.20) $$
we have 
$$ a_n \le C. \eqno (2.21) $$

\medskip

\noindent Thus it follows from (2.19) and (2.21) that 

$$\int_{B_{R}(0)}e^{q \widehat{u}_n } dV_n \le C,  $$
which gives 
$$\int_{B_{R}(0)}e^{q \widehat{u}_n } dx \le C \eqno (2.22) $$
since $g_n$ is bounded in $C^k(B_R(0)) \ \forall \ k  \in \mathbb{N}$, where $dx$ is the Euclidean volume element of $\mathbb{R}^4$. 

\medskip

\medskip

Now, we have by (2.14) since $R \ge 1$ 
$$\rho = \int_{B_{r_n}(x_n)}e^{4u_n} dV_0 = \int_{B_1(0)}e^{4\widehat{u}_n} dV_n \le \int_{B_R(0)}e^{4\widehat{u}_n} dV_n= e^{a_n}\int_{B_R(0)} e^{4(\widehat{u}_n -a_n)} dV_n, $$
and since by (2.19) and H\"older's inequality  we have $\displaystyle \int_{B_R(0)} e^{4(\widehat{u}_n -a_n)} dV_n  \le C$ (recall here that $g_n$ is bounded in $C^k(B_R(0)) \ \forall \ k \in \mathbb{N}$), then we obtain 
$$e^{a_n} \ge  C^{-1}\rho  , $$
which together with (2.21) give 
$$|a_n| \le C. \eqno (2.23)$$

\medskip

It follows from (2.15) and (2.23),  that $(\widehat{u}_n)_n$ is bounded in  $W^{3, p}(B_R(0))$ (for all $p \in [1,  4/3)$). Then by  using  Rellich-Kondrachov  Theorem  and  passing to a subsequence,   we have that 
$(\widehat{u}_n)_n$ converges strongly  in $W^{2, \alpha}(B_R(0))$ for all $\alpha \in [1,2)$,   and in  $L^{\beta}(B_R(0))$ for all $\beta \in [1, + \infty)$,  to a function $\widehat{u}_{\infty} \in W^{3, p}(B_R(0))$ for all $p\in [1, 3/4)$.  That is, 

$$ \forall \ \alpha \in [1, 2), \ \   \|\widehat{u}_n - \widehat{u}_{\infty} \|_{W^{2,\alpha}(B_R(0))}   \underset{n \to + \infty}{\longrightarrow} 0   \eqno (2.24) $$
and 
$$ \forall \beta \in [1, + \infty), \   \  \|\widehat{u}_n - \widehat{u}_{\infty} \|_{L^{\beta}(B_R(0))}  \underset{n \to + \infty}{\longrightarrow} 0 . \eqno (2.25)  $$

\medskip

\noindent Thus it follows from  (2.22) and (2.25) by using H\"older inequality (recall that $q >4$ in (2.22))  that 
$$\int_{B_{R}(0)}\left| e^{4 \widehat{u}_n } - e^{ 4\widehat{u}_{\infty}}\right|  dx   \underset{n \to + \infty}{\longrightarrow} 0  . \eqno (2.26)$$

\medskip

Since $\widehat{u}_n$ satisfies equation (2.9), then by passing to the limit in this equation (in the distributional sense) where we use  (2.26), we obtain 
$$\Delta^2 \widehat{u}_{\infty} = k_0 e^{4\widehat{u}_{\infty}}   \  \  \hbox{in} \  \mathbb{R}^4 ,  \eqno (2.27) $$
where $\Delta$ is the Laplacian in $\mathbb{R}^4$ with respect to the Euclidean metric. By using (2.20), we see that $\widehat{u}_{\infty}$ satisfies also 
$$\int_{B_R(0)}e^{4\widehat{u}_{\infty}(z)} dz \le 1, $$
and since $R \ge 1$ is arbitrary, then 
$$\int_{\mathbb{R}^4}e^{4\widehat{u}_{\infty}(z)} dz \le 1.   \eqno (2.28)$$
\medskip

The solutions of equation (2.27)  satisfying (2.28) are classified in \cite{cL}.  More precisely, it is proved in \cite{cL} that either 
$$\widehat{u}_{\infty}(z) = \log\left( {2\lambda  \over 1 +  \lambda^2  |z-z_0|^2 }\right) - {1\over 4}\log\left({k_0 \over 6}\right)  \eqno (2.29) $$
for some $z_0 \in \mathbb{R}^4$ and $\lambda >0$, or one has
$$- \Delta\widehat{u}_{\infty}(z)  \underset{|z| \to + \infty}{\longrightarrow} a \eqno (2.30)  $$ 
for some  $a >0$. But  by using Proposition 2.1 and (2.24) (by taking $\alpha=1$), one can easily check that  for any $R \ge 1$,
$$\int_{B_R(0)}|\Delta \widehat{u}_{\infty}(z)| dz \le CR^2.  \eqno(2.31)$$
On the other hand,  if (2.30) occurs, then one has for $R $ large enough 
$$ \int_{B_R(0)}|\Delta \widehat{u}_{\infty}(z)|  dz \ge C a R^4$$
which contradicts (2.31). This proves that $ \widehat{u}_{\infty} $  is of the form (2.29). The proof of Proposition 2.4 is then complete.

\end{proof}

\medskip

We close this section with the following well known Adam's inequality (see \cite{aC}) :

\medskip

\begin{prop} Let $(M,g_0)$ be a compact Riemannian $4$-manifold whose  Paneitz  operator $P_0$ is positive  with trivial kernel. Then for any $u \in H^2(M)$,  we have 
$$\int_M e^{4(u - \overline{u}) } dV_0 \le C \exp \left({1\over 4\pi^2} \int_M P_0u\cdot u \ dV_0 \right), \eqno (2.32) $$
where $C$ is a positive constant depending only on $(M, g_0)$ and where \  $\displaystyle \overline{u} = {1\over |M|}\int_M u \ dV_0$ is the average of $u$ on $M$. 
\end{prop}

\bigskip

\section{Integral  Harnack type  inequality }

\medskip

In this section we shall prove an integral  Harnack type inequality, which  is an important tool in the proof of our results. In what follows, we set  $f^+ = \max(f, 0)$ and 
$f^- =  \max(- f, 0)$  for  any function $f$ on $M$.  As already noticed in Section 2, in order to simplify the notations, we are denoting by $|x-y|$ the Riemannian distance between two points $x, y \in M$.  The diameter of $M$ is denoted by $diam(M)$. 

\medskip

 \begin{prop}  Let $ h \in C^0(M)$, and let $(u_n , h_n)_n$ be a sequence in $H^2(M)\times L^1(M)$ satisfying 
 
 \medskip
 
 $$ P_0 u_n + h = h_n , \eqno (3.1) $$
 where we suppose that  
 
 $$\lim_{n\to + \infty} \int_M h_n^{-} dV_0 = 0.    \eqno (3.2)$$
 
 \medskip
 
 \noindent Let  $(x_n , y_n)_n $ a sequence  in $M\times M$,  and let  $(R_n)_n$ with  $0<  R_n  \le \hbox{diam}(M)$  a sequence of positive numbers    satisfying,  for some constant $C_0$ independent of $n$, 
 
 $$|x_n - y_n|  \le C_0 R_n \   \   \     and  \  \   \   \int_{B_{2R_n}(y_n) } h_n^{+} dV_0 \le  \pi^2. \eqno (3.3) $$
Then  for any sequence $(r_n)_n$  such that $0 < r_n \le R_n$ we have 
 $$\int_{B_{R_n}(y_n)}e^{4u_n} dV_0 \le C \left({r_n\over R_n}\right)^{-4 + {1\over 2\pi^2}\|h_n^{+}\|_{L^1(B_{r_n}(x_n))}  + o_n(1) }\int_{B_{r_n}(x_n)}e^{4u_n} dV_0 ,  \eqno (3.4) $$
 where $o_n(1)  \to  0 $ as $n \to + \infty$, and where $C$ is a positive constant independent of $n$. 
 
 \end{prop}

\bigskip

\begin{rem}  As it can be seen in the proof, one can replace $\pi^2$ in   (3.3) in the above proposition   by any positive  constant $\rho < 4 \pi^2$. 

\end{rem} 

\medskip

   \begin{proof}   In what follows,  $C$ is a positive constant independent of $n$  whose values may change from  line to line. Also, to simplify the notations,  we set $R= R_n$ and $r= r_n$. 
   First, let us recall the asymptotic formula for the Green function (see section 2) 
   $$G(x, y) = -  {1 \over 8 \pi^2} \log|x-y|  +  O(1) , \  \  \text{for all} \  x \not= y \  \text{in} \  M .  \eqno (3. 5)  $$
   
   \noindent From the Green representation formula (see section 2)  we have, for any $x \in M$, 
   
   $$u_n(x)- {\bar u}_n = \int_M G(x, y)h_n(y) dV_0(y)- \int_M G(x, y)h(y) dV_0(y)$$
   $$=  \int_M G(x, y)h_n^{+}(y) dV_0(y)- \int_M G(x, y)h_n^{-}(y) dV_0(y)- \int_M G(x, y)h(y) dV_0(y).  \eqno (3.6)$$
   
   \medskip
   
   \noindent Since $h \in C^0(M)$,  then by using (3.5) we have 
   
   $$\left|\int_M G(x,y)h(y) dV_0(y) \right| \le C \|h\|_{L^{\infty}(M)} . \eqno (3.7) $$
   
   \medskip
   
   \noindent  Thus  it follows from (3.6) and (3.7) that for any $x \in M$, 
   $$u_n(x) - {\bar u}_n   \le C \|h\|_{L^{\infty}(M)} + \int_{B_{2R}(y_n)}G(x,y) h_n^+(y)dV_0(y) 
   +  \int_{M\setminus B_{2R}(y_n)}G(x,y)h_n^+(y)dV_0(y)$$
   $$ -  \  \int_M G(x,y)h_n^{-}(y)dV_0(y)$$

   \medskip
   
   \medskip

  \noindent  which implies by integrating the function $e^{4(u_n(x)- {\bar u}_n)}$ on $B_R(y_n)$, 
  
  \medskip
  
  $$\int_{B_R(y_n)}e^{4(u_n - {\bar u}_n)} dV_0 \le e^{C \|h\|_{L^{\infty}(M)}}  \int_{B_R(y_n)} \exp\left(\int_{B_{2R}(y_n)}4G(x,y) h_n^+(y)dV_0(y) \right) dV_0(x) $$
   $$ \times \   \exp\left( \sup_{x\in B_R(y_n)}\int_{M\setminus B_{2R}(y_n)}4 G(x,y)h_n^+(y)dV_0(y) - 
   \inf_{z\in B_R(y_n)}\int_M 4G(z,y)h_n^{-}(y)dV_0(y)\right). \eqno (3.8)  $$

  \bigskip

   \noindent By Jensen inequality, we have 
   
   $$  \exp\left(\int_{B_{2R}(y_n)}4G(x,y) h_n^+(y)dV_0(y)\right)  $$
   $$ \le \  {1 \over  \|h_n^+\|_{L^1(B_{2R}(y_n))}} \int_{B_{2R}(y_n)}  h_n^+(y) \exp\Bigl(4G(x,y)\|h_n^+\|_{L^1(B_{2R}(y_n))}\Bigr) dV_0(y),  $$
   
   \medskip
   
   \medskip
   
\noindent and  integrating this inequality on $B_R(y_n)$ (in the $x$-variable) and using Fubini's Theorem,  we obtain 

$$\int_{B_R(y_n)} \exp\left(\int_{B_{2R}(y_n)}4G(x,y) h_n^+(y)dV_0(y) \right) dV_0(x) $$
$$\le \  {1 \over  \|h_n^+\|_{L^1(B_{2R}(y_n))}} \int_{B_{2R}(y_n)} \left(\int_{B_R(y_n)} \exp\Bigl(4G(x,y)\|h_n^+\|_{L^1(B_{2R}(y_n))}\Bigr) dV_0(x) 
\right)  h_n^+(y)  dV_0(y)   $$
$$\le \ \sup_{y \in B_{2R}(y_n)} \int_{B_R(y_n)} \exp\Bigl(4G(x,y)\|h_n^+\|_{L^1(B_{2R}(y_n))}\Bigr) dV_0(x) . \eqno (3.9)  $$

\medskip

\noindent It follows  from (3.9)   and (3.5)  that 

$$\int_{B_R(y_n)} \exp\left(\int_{B_{2R}(y_n)}4G(x,y) h_n^+(y)dV_0(y) \right) dV_0(x) $$ 
$$\le \    \exp\Bigl( C \|h_n^+\|_{L^1(B_{2R}(y_n))}\Bigr)   \sup_{y \in B_{2R}(y_0)} \int_{B_R(y_n)} |x-y|^{- {1\over 2\pi^2} \|h_n^+\|_{L^1(B_{2R}(y_n))}} dV_0(x)$$
$$ \le e^{\pi^2C}\sup_{y \in B_{2R}(y_n)} \int_{B_R(y_n)} |x-y|^{- {1\over 2\pi^2} \|h_n^+\|_{L^1(B_{2R}(y_n))}} dV_0(x), \eqno (3.10) $$

\medskip

\noindent where we have used    (3.3). 

\medskip

\medskip

\noindent It is easy to check that,  for any $ \alpha \in [0,  4)   $, and any $y \in M$,  we have 
$$ \int_{B_R(y_n)} |x-y|^{ - \alpha} dV_0(x)  \le {C \over 4-\alpha} R^{4-\alpha  }.  $$

\noindent  Since $  \|h_n^+\|_{L^1(B_{2R}(y_n))} \le  \pi^2$ by (3.3), it follows by taking $\alpha =  {1\over 2\pi^2} \|h_n^+\|_{L^1(B_{2R}(y_n))} \le 1$, 

$$ \int_{B_R(y_n)} |x-y|^{- {1\over 2\pi^2} \|h_n^+\|_{L^1(B_R(y_n))}} dV_0(x)  \le C R^{4- {1\over 2\pi^2} \|h_n^+\|_{L^1(B_{2R}(y_n))}} 
 .  \eqno (3.11)$$

\noindent Thus it follows from (3.8), (3.10)  and (3.11) that 

$$\int_{B_R(y_n)}e^{4(u_n - {\bar u}_n)} dV_0 \le C  R^4 \exp\left(   C \|h\|_{L^{\infty}(M)}  - { \log R \over 2\pi^2} \|h_n^+\|_{L^1(B_{2R}(y_n))}\right) $$
$$ \times  \exp\left( \sup_{x\in B_R(y_n)}\int_{M\setminus B_{2R}(y_n)}4 G(x,y)h_n^+(y)dV_0(y) - 
   \inf_{z\in B_R(y_n)}\int_M 4G(z,y)h_n^{-}(y)dV_0(y)\right). \eqno (3.12) $$
   
   \bigskip
   
   \medskip

   On the other hand, using again   the representation formula (3.6), we have by using (3.7), for any $x \in B_r(x_n)$, 
   
   $$u(x)- {\bar u}_n = \int_M G(x, y)h_n^+(y) dV_0(y)-  \int_M G(x, y)h_n^-(y) dV_0(y) - \int_M G(x, y)h(y) dV_0(y)$$
   $$\ge  \int_{B_r(x_n)} G(x, y)h_n^+(y) dV_0(y) +  \int_{M\setminus B_r(x_n)} G(x, y)h_n^+(y) dV_0(y) - \int_M G(x, y)h_n^-(y) dV_0(y)  - C\|h\|_{L^{\infty}(M)}$$
$$ \ge \inf_{z \in B_r(x_n)}\int_{B_r(x_n)} G(z, y)h_n^+(y) dV_0(y) + \inf_{z \in B_r(x_n)}\int_{M\setminus B_r(x_n)} G(z, y)h_n^+(y) dV_0(y) $$
$$ - \  \int_M G(x, y)h_n^-(y) dV_0(y)  -  C\|h\|_{L^{\infty}(M)}. \eqno (3.13)$$

\noindent Since by  (3.5) we have $G(z,y) \ge - {1\over 8\pi^2}\log r - C$ for any $z , y \in B_r(x_n)$, then it follows from (3.13) that for any $x \in B_r(x_n)$, 

$$e^{4(u_n(x)-{\bar u}_n )} \ge  \exp\left(- C\|h\|_{L^{\infty}(M)} -C\|h_n^+\|_{L^1(M)}   - {\log r \over 2\pi^2} \|h_n^+\|_{L^1(B_r(x_n))}\right) $$
$$ \times  \exp\left( \inf_{z \in B_r(x_n)}\int_{M\setminus B_r(x_n)} 4 G(z, y)h_n^+(y) dV_0(y)\right)   \exp\left( - \int_M 4G(x, y)h_n^-(y) dV_0(y)  \right).  \eqno (3.14) $$

\noindent But from (3.1) and (3.2) we have $\|h_n^+\|_{L^1(M)} \le  \|h\|_{L^1(M)}  + o_n(1) \le C \|h\|_{L^{\infty}(M)} + o_n(1)$.  So it follows from (3.14)  on  integrating  on $B_r(x_n)$ that  

$$\int_{B_r(x_n)} e^{4(u_n - {\bar u}_n )}dV_0  \ge C \exp\left( - C\|h\|_{L^{\infty}(M)}  - {\log r \over 2\pi^2} \|h_n^+\|_{L^1(B_r(x_n))}\right)  $$
$$ \times  \exp\left( \inf_{z \in B_r(x_n)}\int_{M\setminus B_r(x_n)} 4 G(z, y)h_n^+(y) dV_0(y)\right)  \int_{B_r(x_n) }\exp\left( - \int_M 4 G(x, y)h_n^-(y) dV_0(y)  \right) dV_0(x).  \eqno (3.15) $$

\medskip

\noindent But by Jensen inequality we have 
$$ \int_{B_r(x_n)} \exp\left( - \int_M 4 G(x, y)h_n^-(y) dV_0(y)  \right) dV_0(x) $$
$$ \ge |B_r(x_n)| \exp\left( - {1 \over |B_r(x_n)|} \int_{B_r(x_n)} \int_M 4 G(x,y)
h_n^-(y) dV_0(y) dV_0(x)  \right) . \eqno (3.16) $$

\medskip

\noindent Since $ |B_r(x_n)| \ge Cr^4$, it follows from (3.15) and (3.16) that 

$$\int_{B_r(x_n)} e^{4(u_n - {\bar u}_n )}dV_0  \ge C r^4 \exp\left(  - C\|h\|_{L^{\infty}(M)} - {\log r \over 2\pi^2} \|h_n^+\|_{L^1(B_r(x_n))}\right)$$
$$  \times  \exp\left( \inf_{z \in B_r(x_n)}\int_{M\setminus B_r(x_n)} 4 G(z, y)h_n^+(y) dV_0(y)
  - {1 \over |B_r(x_n)|} \int_{B_r(x_n)} \int_M 4G(x,y) h_n^-(y) dV_0(y) dV_0(x) \right).   \eqno (3.17)$$

\medskip

Now since 
$$   { \int_{B_R(y_n)}e^{4u_n } dV_0 \over  \int_{B_r(x_n)} e^{4 u_n }dV_0}  =  { \int_{B_R(y_n)}e^{4(u_n - {\bar u}_n)} dV_0 \over  \int_{B_r(x_n)} e^{4(u_n - {\bar u}_n )}dV_0} , $$
then it follows  from (3.12) and (3.17) that 
$$   { \int_{B_R(y_n)}e^{4u_n } dV_0 \over  \int_{B_r(x_n)} e^{4 u_n }dV_0}   \le  \  C \left({r\over R}\right)^{- 4} \exp\left( C\|h\|_{L^{\infty}(M)} + {\log r \over 2\pi^2} \|h_n^+\|_{L^1(B_r(x_n))} -  {\log R \over 2\pi^2} \|h_n^+\|_{L^1(B_{2R}(y_n))}\right)
$$
   $$\times \exp\left( \sup_{x\in B_R(y_n)}\int_{M\setminus B_{2R}(y_n)}4 G(x,y)h_n^+(y)dV_0(y) - 
   \inf_{z\in B_r(x_n)}\int_{M\setminus B_r(x_n)}  4G(z,y)h_n^{+}(y)dV_0(y)\right)$$
   $$\times  \exp\left(  {1 \over |B_r(x_n)|} \int_{B_r(x_n)} \int_M 4G(x,y)
h_n^-(y) dV_0(y) dV_0(x)  -    \inf_{z\in B_R(y_n)}\int_M 4G(z,y)h_n^{-}(y)dV_0(y) \right)  . \eqno (3.18) $$

\medskip

\medskip

Set 
$$A =  \exp\left( \sup_{x\in B_R(y_n)}\int_{M\setminus B_{2R}(y_n)}4 G(x,y)h_n^+(y)dV_0(y) - 
   \inf_{z\in B_r(x_n)}\int_{M\setminus B_r(x_n)}  4G(z,y)h_n^{+}(y)dV_0(y)\right)$$
   and 
   $$ B = \exp\left(  {1 \over |B_r(x_n)|} \int_{B_r(x_n)} \int_M 4G(x,y)
h_n^-(y) dV_0(y) dV_0(x)  -    \inf_{z\in B_R(y_n)}\int_M 4G(z,y)h_n^{-}(y)dV_0(y) \right) .  $$

\medskip

\noindent We shall prove that

$$A \le C \exp\left( {\log R \over 2\pi^2} \|h_n^+\|_{L^1(B_{2R}(y_n))} -  {\log R \over 2\pi^2} \|h_n^+\|_{L^1(B_r(x_n))}  \right)  \eqno (3.19)$$

\noindent and

$$B \le C  \exp\left(   C \|h_n^-\|_{L^1(M)}\log {R\over r }   \right) . \eqno (3.20) $$

\medskip

\noindent  It is clear that Proposition 3.1 will follow from (3.18), (3.19) and (3.20) by using (3.2). Let us then prove the estimates (3.19) and (3.20). 

\medskip

First we shall prove (3.19). We have for any $x \in B_R(y_n)$ and $z \in B_r(x_n)$
$$\int_{M\setminus B_{2R}(y_n)}4 G(x,y)h_n^+(y)dV_0(y) - 
   \int_{M\setminus B_r(x_n)}  4G(z,y)h_n^{+}(y)dV_0(y) $$
   $$ =  \int_{M\setminus B_{4C_0R}(y_n)} 4 (G(x,y) -G(z,y) ) h_n^+(y)dV_0(y) $$
  $$ + \int_{B_{4C_0 R}(y_n)\setminus B_{2R}(y_n)} 4G(x,y) h_n^+(y)dV_0(y)  -  \int_{B_{4C_0R}(y_n)\setminus B_r(x_n)}  4G(z,y)h_n^{+}(y)dV_0(y)  \eqno (3.21) $$
   
   \medskip
   
  \noindent  since  $B_{2R}(y_n) \subset  B_{4C_0 R}(y_n)$ and  $B_{r}(x_n) \subset  B_{4C_0 R}(y_n)$, where $C_0$ is the constant in (3.3) that we assume satisfying $C_0 \ge 1$ without loss of generality  (we recall here that  $r = r_n \le R_n= R$). 
  
  \medskip 
  
  Let us estimate the first  term in the right side of  (3.21). We have  for any $x \in B_R(y_n)$, $z \in B_r(x_n) $ and $y \in M\setminus B_{4 C_0 R}(y_n)$ by using  the  hypothesis  $|x_n - y_n| \le C_0 R$ and $r \le R$, that 
 
  $$|x-y| \ge \ |y-y_n| - |x-y_n|\ \ge 4C_0R - R  \ge  3 C_0 R  ,  \eqno (3.22) $$
  
  \noindent and

 $$|z-y| \ge \ |y-y_n| - |y_n - x_n|- |z-x_n| \ \ge 4C_0 R - C_0 R - r   \ge  2 C_0 R .   \eqno (3.23) $$
 
 \medskip
 
On the other hand, we have by using (3.22) 
 
 $$|z-y|   \le |z-x_n| + |x_n-y_n| + |y_n- x|  + |x-y| $$
 $$\le r + C_0 R + R + |x-y|  \le 3 C_0 R  + |x-y| $$
 $$ \le  \  2  |x-y|,  $$
 and by using (3.23) we have 
 $$ |x-y|  \le |x-y_n| + |x_n- y_n| + |x_n - z |  + |z-y|  $$
 $$\le \ R + C_0 R + r   + |z-y| \le  3 C_0 R  + |z-y|   $$
 $$\le \ {5\over 2} |z-y| . $$
 
 \noindent Thus we have 
 $${2\over 5} \le  {|z-y| \over |x-y|} \le 2  . \eqno (3.24) $$
 
 \noindent It follows from (3.5) and (3.24) that for any $x \in B_R(y_n)$, $z \in B_r(x_n) $ and $y \in M\setminus B_{4 C_0 R}(y_n)$, 
 $$\left|G(x, y)-G(z, y)\right| \le C , \eqno (3.25) $$
 
 \noindent which gives 
 $$\int_{M\setminus B_{4C_0 R}(y_n)}  4\Bigl( G(x,y) -G(z,y)\Bigr) 
h_n^+(y) dV_0(y)   \le C \|h_n^+\|_{L^1(M)}   \eqno (3.26)$$

\medskip

\noindent for any $x \in B_R(y_n)$, $z \in B_r(x_n) $.

\medskip

 Now we shall estimate the second and third  term in the right side (3.21).  By using  formula (3.5) we have for any $x \in B_R(y_n)$  and 
 $y  \in B_{4C_0 R}(y_n)\setminus B_{2R}(y_n)$, 
 $$G(x, y ) \le -  { 1\over 8\pi^2}\log R + C .  \eqno (3.27)$$
 We have for any $z \in B_r(x_n)$  and  $y  \in B_{4C_0 R}(y_n)\setminus B_{r}(x_n)$, since $|x_n-y_n|  \le C_0 R$ (by hypothesis)  and $r \le R$, 
 $$|y-z| \le  |y-y_n| + |y_n - x_n| + |x_n- z| \le 4 C_0R + C_0 R + r \le 6 C_0 R $$
 which implies by (3.5) that 
 $$G(z,y) \ge-  {1 \over 8\pi^2}\log R + C .  \eqno (3.28)$$

 \medskip
 
 \noindent  It follows from (3.27) and (3.28),  for any $x \in B_R(y_n)$ and $z \in B_r(x_n)$,  that 
 
   $$ \int_{B_{4C_0 R}(y_n)\setminus B_{2R}(y_n)} 4G(x,y) h_n^+(y)dV_0(y)  -  \int_{B_{4C_0 R}(y_n)\setminus B_r(x_n)}  4G(z,y)h_n^{+}(y)dV_0(y) $$
   $$ \le {\log R \over 2\pi^2}\int_{B_{2R}(y_n)} h_n^+  dV_0 -   {\log R \over 2\pi^2}\int_{B_r(x_n)} h_n^+  dV_0 + C \|h_n^+\|_{L^1(M)}.  \eqno (3.29)$$
   
   \medskip
   
   Combining (3.21), (3.26) and (3.29) we obtain the desired estimate (3.19) since  $ \|h_n^+\|_{L^1(M)} \le \|h\|_{L^1(M)} + \|h_n^-\|_{L^1(M)} \le C$ by integrating (3.1) and using  (3.2). 
   
   \bigskip
   
   Now it remains to prove (3.20). We have 
$$B = \exp\left\{  {1 \over |B_r(x_n)|} \int_{B_r(x_n)} \int_M 4G(x,y)
h_n^-(y) dV_0(y) dV_0(x)  -    \inf_{z\in B_R(y_n)}\int_M 4G(z,y)h_n^{-}(y)dV_0(y) \right\}  $$
$$ =   \sup_{z\in B_R(y_n)}\exp\left\{  {1 \over |B_r(x_n)|} \int_{B_r(x_n)} \left(\int_M 4\Bigl( G(x,y) -G(z,y)\Bigr) 
h_n^-(y) dV_0(y) \right)dV_0(x)  \right\} .  \eqno (3.30)$$

\medskip

\noindent But we have for any $z \in B_R(y_n)$,
$$ {1 \over |B_r(x_n)|} \int_{B_r(x_n)}\left( \int_M 4\Bigl( G(x,y) -G(z,y)\Bigr) 
h_n^-(y) dV_0(y) \right) dV_0(x) $$
$$ =  {1 \over |B_r(x_n)|} \int_{B_r(x_n)}\left( \int_{B_{4C_0 R}(y_n)} 4\Bigl( G(x,y) -G(z,y)\Bigr) 
h_n^-(y) dV_0(y)\right) dV_0(x)  $$
$$ +  \   {1 \over |B_r(x_n)|}\int_{B_r(x_n)}\left( \int_{M\setminus B_{4C_0 R}(y_n)}  4\Bigl( G(x,y) -G(z,y)\Bigr) 
h_n^-(y) dV_0(y) \right) dV_0(x) . \eqno (3.31) $$

\medskip

\noindent Since  by (3.5) we have,   for any $z \in B_R(y_n)$  and  $y \in B_{4C_0 R}(y_n)$,  
 $$G(z,y)  \ge -{1 \over 8\pi^2}\log R + C,  $$ 
   then the first term in the right side of  (3.31) can be estimated as follows 
$${1 \over |B_r(x_n)|} \int_{B_r(x_n)} \int_{ B_{4C_0 R}(y_n)} 4\Bigl( G(x,y) -G(z,y)\Bigr) 
h_n^-(y) dV_0(y) dV_0(x) $$
$$ \le  \  {1 \over |B_r(x_n)|}  \int_{B_r(x_n)} \int_{  B_{4C_0 R}(y_n)} {1 \over 8\pi^2}\Bigl( -  \log|x-y|  +  \log R \Bigr) 
h_n^-(y) dV_0(y) dV_0(x)  + C \|h_n^{-}\|_{L^1(M)}  $$
  $$\le{1 \over 8 \pi^2}\|h^{-}_n\|_{L^1(M)}   {1 \over |B_r(x_n)|} \sup_{y\in  B_{4C_0 R}(x_n)} \int_{B_r(x_n)} \left|\log\left({1\over R}|x-y|\right)\right| dV_0(x) + C \|h_n^{-}\|_{L^1(M)} .  \eqno (3.32)$$

 \medskip
 
 \noindent A direct computation shows that 
 
 $$  {1 \over |B_r(x_n)|} \sup_{y\in B_{4C_0R}(y_n)} \int_{B_r(x_n)} \left|\log\left({1\over R}|x-y|\right)\right| dV_0(x) \le C\log\left({R\over r} \right) + C  . $$

\noindent (we recall here that $r \le R$).  Combining the last inequality  with (3.32) gives for any $z \in B_R(y_n)$, 

 $$ {1 \over |B_r(x_n)|}\int_{B_r(x_n)} \int_{B_{4C_0R}(y_n)} 4\Bigl( G(x,y) -G(z,y)\Bigr) 
h_n^-(y) dV_0(y) dV_0(x) $$
$$\le  C \|h^{-}_n\|_{L^1(M)} \log\left({R\over r} \right) + C  \|h^{-}_n\|_{L^1(M)}  . \eqno (3.33) $$

 Now let us estimate  the second term  in the right side of (3.31).  We recall that from  (3.25)  we have  
 $$ \left|G(x, y)-G(z, y)\right| \le C ,  $$
  for any  $x \in B_r(x_n), z \in B_R(y_n)$  and $y \in M\setminus B_{4 C_0R}(y_n)$. Thus we obtain 

 $${1 \over |B_r(x_n)|}\int_{B_r(x_n)} \int_{M\setminus B_{4C_0 R}(x_0)}  4\Bigl( G(x,y) -G(z,y)\Bigr) 
h_n^-(y) dV_0(y) dV_0(x)  \le C \|h_n^-\|_{L^1(M)}   . \eqno (3.34)$$

\noindent It follows from  (3.30), (3.31), (3.33) and (3.34) that 
$$  B \le \exp\left( C \|h^{-}_n\|_{L^1(M)} \log\left({R\over r} \right) + C \|h_n^-\|_{L^1(M)}\right)  , $$
which proves (3.20) by using (3.2). The proof of  Proposition 3.1 is then complete.

\end{proof}

\bigskip

\section{Volume quantization and proof of Theorem 1.1 }

\medskip

In this section we apply the result of section 3 (Harnack type inequality) to prove some fundamental properties on solutions of equation (1.8) in Theorem 1.1.  They state that  the conformal volume concentrates with quantization  at  points corresponding to blow-up sequences.  Through the rest of the paper we shall assume that $k_0 = \int_M f dV_0  > 0$ where $ f$ is as in Theorem 1.1. Indeed,  if $k_0 \le 0$, then Theorem 1.1  is a direct consequence of Proposition 2.3 in section 2.  

\medskip

\begin{prop}  Let $(u_n)_n$ as in Theorem 1.1 and let $(x_n, r_n)_n$ a blow-up for the sequence  $(u_n)_n$. Let $(y_n)_n$ a sequence of points in $M$, and $ 0< \rho_n \le  \hbox{diam}(M)$ such that \  $\displaystyle \lim_{n\to+\infty}{r_n\over \rho_n} = 0$. Suppose that, for some positive constant $C_0$ independent of $n$, we have 
$$ |x_n-y_n| \le C_0 \rho_n   \   \    \hbox{and} \  \    \int_{B_{2\rho_n}(y_n)} e^{4u_n} dV_0 \le { \pi^2 \over k_0}  \ , $$
where  $\displaystyle k_0 = \int_M f dV_0$.  Then 
$$  \int_{B_{\rho_n}(y_n)} e^{4u_n} dV_0 \le C \left(r_n \over \rho_n\right)^{2 + o_n(1)} , \eqno (4.1)$$
where $C$ is a positive constant independent of $n$. In particular, we have 

$$\lim_{n\to + \infty}  \int_{B_{\rho_n}(y_n)} e^{4u_n} dV_0  =  0.$$

\end{prop} 

\bigskip

\begin{proof}

First let us apply Proposition 2.4  by choosing $\beta_n= \sqrt{\rho_n\over r_n}$. Then there exists $b_n  \le \sqrt{\rho_n\over r_n}$ such that $b_n \to + \infty$ and 
$$\lim_{n\to +\infty} \int_{B_{b_nr_n(x_n)}} e^{4u_n} dV_0 = {16\pi^2 \over k_0}.  \eqno (4.2)$$

\medskip

We can apply  now Proposition 3.1 to $(u_n)_n$ by choosing $h_n= k_ne^{4u_n}-f_n+ f$, $ h = f$,   $R_n =  \rho_n$, and $b_nr_n$ instead of $r_n$. Indeed,  since $f_n \underset{n \to + \infty}{\longrightarrow}  f$ in $L^1(M)$ and $k_0 = \int_M fdV_0 > 0$, then  $k_n= \int_Mf_n dV_0  > 0 $  for $n$ large enough. Then  one can  easily check that hypothesis (3.2)-(3.3) in Proposition 3.1  are  satisfied. Thus we obtain
 $$ \int_{B_{\rho_n}(y_n)} e^{4u_n} dV_0 \le C\left({b_nr_n\over \rho_n}\right)^{-4 + {1\over 2\pi^2}\|h_n^{+}\|_{L^1(B_{b_nr_n}(x_n))}  + o_n(1) }\int_{B_{b_nr_n}(x_n)}e^{4u_n} dV_0$$
 $$\le C\left({b_nr_n\over \rho_n}\right)^{-4 + {1\over 2\pi^2}\|h_n^{+}\|_{L^1(B_{b_nr_n}(x_n))}  + o_n(1) }, \eqno (4.3) $$
 where we have used the fact that $ \displaystyle  \int_{B_{b_nr_n}(x_n)}e^{4u_n} dV_0 \le \int_M e^{4u_n} dV_0 =  1.$ 
 
 Since $h_n  = k_ne^{4u_n}-f_n+f $, then we have by using (4.2)  and the fact that $f_n \to f$ in $L^1(M)$, that 
 $$ {1\over 2\pi^2}\|h_n^{+}\|_{L^1(B_{b_nr_n}(x_n))}= 8 + o_n(1) $$
 and by replacing in  (4.3)we get 
 $$\int_{B_{\rho_n}(y_n)} e^{4u_n} dV_0 \le   C\left({b_nr_n\over \rho_n}\right)^{4+ o_n(1) }. $$
  This proves  estimate (4.1) since  $b_n \le \sqrt{\rho_n \over r_n}$.

 \end{proof}
 
 \bigskip

\begin{prop} Let $(u_n)_n$ as in Theorem 1.1 and let $(x_n, r_n)_n$ a blow-up for the sequence  $(u_n)_n$. Let    \ $0 < R_n \le S_n$ such  that ${r_n \over R_n}  \underset{n \to + \infty}{\longrightarrow}  0 $, and suppose that    there exists  a positive constant  \ $\alpha \le 1 $ independent of $n$  such that
$$ \forall \ B_r(y) \subset  B_{2S_n}(x_n)\setminus B_{{1\over 2}R_n}(x_n), \   \int_{B_r(y)} e^{4u_n} dV_0 \ge { \pi^2\over k_0}  \  \Longrightarrow \   r \ge \alpha \hspace{0,5mm} |y-x_n|  . \eqno (4.4) $$
Then 
$$\lim_{n \to + \infty} \int_{B_{S_n}(x_n)\setminus B_{R_n}(x_n)}e^{4u_n} dV_0 =  0. \eqno (4.5) $$
\end{prop} 
 
 \medskip
 
 \begin{proof}  Before giving the proof we note here that we may assume without loss of generality that $S_n \le \hbox{diam}(M)$. First we shall prove  that for any \ $ \rho_n \in [ R_n , S_n]$ we have the following estimate 
 $$\int_{B_{{3\over 2}\rho_n}(x_n)\setminus B_{\rho_n}(x_n)}e^{4u_n} dV_0 \le C\left({r_n\over \rho_n}\right)^{2+ o_n(1)}, \eqno (4.6) $$
 where $C$ is constant independent of $n$, and $o_n(1) \to 0$ as $n \to + \infty$.  Then  (4.5) will follow from (4.6) by using  an appropriate decomposition of the annulus $B_{S_n}(x_n)\setminus B_{R_n}(x_n)$. Indeed, suppose that (4.6) is proved, then by choosing $N \in \mathbb{N}$ such that $ (3/2)^N \le S_n/R_n \le (3/2)^{N+1}$,  and applying (4.6) with $\rho_n = (3/2)^jR_n$ for $j = 0, ..., N$, we obtain for $n$ large enough :
 
 $$\int_{B_{(3/2)^{j+1}R_n}(x_n)\setminus B_{(3/2)^jR_n}(x_n)}e^{4u_n} dV_0 \le C (2/3)^{(2+o_n(1))j} \left({r_n\over R_n}\right)^{2+ o_n(1)} \le C(2/3)^{j} \left({r_n\over R_n}\right)^{2+ o_n(1)}, $$
 
 \medskip

 \noindent and by summing up  over $j= 0, ..., N$, one gets 
 $$ \int_{B_{S_n}(x_n)\setminus B_{R_n}(x_n)}e^{4u_n} dV_0  \le \sum_{j= 0}^{N} \int_{B_{(3/2)^{j+1}R_n}(x_n)\setminus B_{(3/2)^jR_n}(x_n)}e^{4u_n} dV_0 \le $$
 $$C  \left({r_n\over R_n}\right)^{2+ o_n(1)} \sum_{j= 0}^{N} (2/3)^{j} \le 3 C  \left({r_n\over R_n}\right)^{2+ o_n(1)}  \to 0 \hspace{2mm} \hbox{as} \  \  n \to + \infty, $$
 
 \medskip
 
\noindent  which proves the desired result (4.5).  

\medskip

Now let us prove the estimate (4.6).   We can cover the set      $ B_{{3\over 2}\rho_n}(x_n)\setminus B_{\rho_n}(x_n)  $ by  a finite number of balls $B_{{1\over 4}\alpha\rho_n}(z_1), ...,  B_{{1\over 4}\alpha\rho_n}(z_L)$, where $L \in \mathbb{N} $  is independent of $n$,  and where $\alpha$ is the constant appearing in (4.4),  such that  
$$B_{{1\over 2}\alpha\rho_n}(z_j) \subset   B_{{2}\rho_n}(x_n)\setminus B_{{1\over2}\rho_n}(x_n) \subset B_{2S_n}(x_n) \setminus B_{{1\over 2}R_n}(x_n) , \  \    j= 1, ..., L. $$

\medskip

\noindent But since $ |z_j-x_n|  \ge \rho_n > {1\over 2}\rho_n$,  we have  from (4.4)  that 
$$ \int_{B_{{1\over 2}\alpha\rho_n}(z_j)} e^{4u_n} dV_0 <  { \pi^2 \over k_0}  \  \   \forall \  j= 1, ..., L.  \eqno(4.7) $$
(we recall here that $0 < \alpha \le 1$.)

 We can now apply Proposition 4.1 by taking  $y_n = z_j$  to get 
$$ \int_{B_{{1\over 4}\alpha\rho_n}(z_j)} e^{4u_n} dV_0 \le C\left({r_n\over \rho_n}\right)^{2+ o_n(1)}  \  \   \forall \  j= 1, ..., L ,  $$
 and the estimate (4.6) follows.   This achieves the proof of Proposition 4.2.

 \end{proof} 
 
 \bigskip
 
 \begin{prop}  Let $(u_n, f_n)$ as in Theorem 1.1. Let $(x_n^1, r_n^1)_n , ...., (x_n^m, r_n^m)_n$  be $m$ blow-ups for $(u_n)_n$, and $R_n^1, ..., R_n^m  > 0 $  such that 
  $$ \lim_{n\to + \infty}{R_n^i }= 0 \ \ , \ \  \lim_{n\to + \infty}{r_n^i \over R_n^i} = 0 \   \   \forall \ i= 1, ..., m ,  \eqno (4.8) $$
  and 
  $$  \lim_{n\to + \infty}{R_n^i \over |x_n^j -x_n^i |} = 0 \   \   \forall \ i\not= j \ \hbox{in} \ \{ 1, ..., m\} \  \ \hbox{if} \  m \ge 2.   \eqno (4.9) $$
Let  \ $ \displaystyle S_n  \ge 4 \max_{i\not=j}|x_n^i-x_n^j|  $ 
and   suppose that there exists a positive constant $\alpha \le 1$ independent of $n$ such that 
 $$ \forall \ B_r(y) \subset  \bigcup_{j=1}^{m}B_{2S_n}(x_n^j)\setminus \bigcup_{j=1}^mB_{{1\over 2}R_n^j}(x_n^j), \   \int_{B_r(y)} e^{4u_n} dV_0 \ge { \pi^2  \over k_0}  \  \Longrightarrow \   r \ge \alpha \hspace{0,5mm} 
 d_n(y) , \eqno (4.10) $$
  where $ \displaystyle d_n(y) = \inf_{1\le j\le m} |y-x_n^j| $.   Then 
$$  \lim_{n \to + \infty} \int_{\bigcup_{j=1}^mB_{S_n}(x_n^j)\setminus\bigcup_{j=1}^mB_{R_n^j}(x_n^j)}e^{4u_n} dV_0 = 0 .  \eqno (4.11) $$

  \end{prop}
  
  \bigskip
  
  \begin{proof} It is clear that it suffices to prove (4.11) for a subsequence of $(u_n)_n$. We proceed by induction on $m$. Suppose $m = 1$, then it follows from  Proposition 4.2,  by taking $x_n= x_n^1$  and  $R_n = R_n^1$, that 
 $$   \lim_{n \to + \infty} \int_{B_{S_n}(x_n^1)\setminus B_{R_n^1}(x_n^1)}e^{4u_n} dV_0 = 0 . $$  
 
 Now  let $m \ge 2$ be an integer and suppose that  (4.11) is true for any $l$ blow-ups  with $l \le m$. We shall prove that this is also the case for any  $(m+1)$ blow-ups.  Let then $(x_n^1, r_n^1)_n , ...., (x_n^{m+1}, r_n^{m+1})_n$  be $(m+1)$ blow-ups for $(u_n)_n$  satisfying  (4.9)-(4.10) for some $R_n^i>0, \  i \in \llbracket  1,  m+1\rrbracket$ and $S_n>0$,   that is 
  $$   \lim_{n\to + \infty}{r_n^i \over R_n^i} = 0 \   \   \forall \ i= 1, ..., m+1 , \    \      \lim_{n\to + \infty}{R_n^i \over |x_n^i -x_n^j |} = 0 \   \   \forall \ i\not= j \ \hbox{in} \ \{ 1, ..., m+1\} , \eqno (4.12)  $$
and 
 $$ \forall \ B_r(y) \subset  \bigcup_{j=1}^{m+1}B_{2S_n}(x_n^j)\setminus \bigcup_{j=1}^{m+1}B_{{1\over 2}R_n^j}(x_n^j), \   \int_{B_r(y)} e^{4u_n} dV_0 \ge { \pi^2\over k_0}   \  \Longrightarrow \   r \ge \alpha \hspace{0,5mm} 
 d_n(y) , \eqno (4.13) $$
   where $ \displaystyle d_n(y) = \inf_{1\le j\le m+1} |y-x_n^j| $.

 \bigskip
 
 Let 
 $$ d_n= \inf \left\{ |x_n^i-x_n^j| \  : \  i, j \in \llbracket 1, m+1\rrbracket,  i \not=j \right\} $$
 and 
 $$ D_n =  \sup \left\{ |x_n^i-x_n^j| \  : \  i, j \in \llbracket 1, m+1\rrbracket,  i \not=j \right\}  .$$

 \bigskip
 
\noindent  By passing to a subsequence if necessary, we distinguish two  cases depending on $d_n$ and $D_n$. That is,  we have  either $D_n \le Cd_n$, where $C$ is a positive constant independent of $n$,  or $\displaystyle \lim_{n\to +\infty} {d_n \over D_n} = 0$. 
 
 \bigskip
 
 \medskip

\noindent  {\bf\underline {First case}} :  \  {\bf $D_n \le Cd_n$, where $C$ is a positive constant independent of $n$. } 

\bigskip

  If we apply Proposition 4.2 by taking  $x_n= x_n^i$ and  $ R_n = 4D_n$ (by using (4.13)),  we have  for any $i = 1, ..., m+1$ 
  $$   \lim_{n \to + \infty} \int_{B_{S_n }(x_n^i)\setminus B_{4D_n}(x_n^i)}e^{4u_n} dV_0 = 0 . $$  
 Thus it remains to prove that 
 $$   \lim_{n \to + \infty} \int_{\bigcup_{j=1}^{m+1}B_{4D_n}(x_n^j)\setminus\bigcup_{j=1}^{m+1}B_{R_n^j}(x_n^j)}e^{4u_n} dV_0 = 0 . \eqno (4.14)$$

 \bigskip

\noindent  We have by (4.12) since $D_n \le Cd_n$ that $\displaystyle \lim_{n\to \infty} {R_n^j \over d_n} = 0, \  j = 1, ..., m+1$. Thus if  we apply  Proposition 4.2 by taking $x_n = x_n^j, \ 
 R_n= R_n^j$ and $S_n= {1\over 4}d_n$  (by using (4.13)),  we obtain  
 $$  \lim_{n \to + \infty} \int_{B_{{1\over 4}d_n}(x_n^j)\setminus B_{R_n^j}(x_n^j)}e^{4u_n} dV_0 = 0 \  \  \forall \ j= 1, ..., m+1.  \eqno (4.15) $$
 
 On the other hand,  since $d_n \le D_n \le Cd_n$, we can cover the set $\bigcup_{j=1}^{m+1}B_{4D_n}(x_n^j)\setminus\bigcup_{j=1}^{m+1}B_{{1\over 4}d_n}(x_n^j)$ by a finite number $N$ (independent of $n$) of balls $B_{{1\over 16}\alpha d_n}(z_n^l) , l = 1, ..., N$, where $ 0 <\alpha \le 1$ is the constant appearing in (4.13),  such that  $B_{{1\over 8}\alpha d_n}(z_n^l) \subset  \bigcup_{j=1}^{m+1}B_{2D_n}(x_n^j) \setminus\bigcup_{j=1}^{m+1}B_{{1\over 8}d_n}(x_n^j) $. Then we can apply Proposition 4.1 by taking  $y_n = z_n^l $, $x_n= x_n^1,  \  r_n= r_n^1,  \ \rho_n= 
 {1\over 16}\alpha d_n$,  and by using  (4.13), we obtain 
  $$  \lim_{n \to + \infty} \int_{B_{{1\over 16}\alpha d_n}(z_n^l)}e^{4u_n} dV_0 = 0 \  \  \forall \ l= 1, ..., N .  \eqno (4.16) $$
  It is clear that (4.14) follows from (4.15) and (4.16).

 \bigskip
 
 \noindent  {\bf\underline {Second  case}} :  \  { \bf $\displaystyle \lim_{n\to +\infty} {d_n \over D_n} = 0 $. } 
 
 \bigskip

   By relabelling the blow-ups and  passing  to a subsequence if necessary, we may suppose that $d_n = |x_n^1-x_n^{2}|$. Define the set $J$ by :
  
 $$J =  \{ \  j \in \llbracket 1, m+ 1\rrbracket  \ : \  |x_n^j-x_n^1| \le C_j d_n   \   \forall n \  \}  \ ,  $$
 
 \bigskip
 
 \noindent where $C_j$  is a positive constant  independent of $n$.  By Taking $\displaystyle C_0= \max_{j\in J} C_j$    we have ( by passing to a subsequence if necessary)
 $$ \forall \ j \in J, \   |x_n^j-x_n^1| \le C_0 d_n  \  \  \forall \ n,  \eqno (4.17) $$
 and 
 $$  \forall j \in 
  \llbracket 1, m+1 \rrbracket \setminus J, \  \lim_{n\to +\infty}{d_n \over |x_n^j -x_n^1 |} = 0  .  \eqno (4. 18)$$

  \medskip

  By relabeling the blow-ups (except for $j= 1$ and $j= 2$ ) and observing that $1, 2  \in J$, we may suppose  that  $J= \llbracket 1, k\rrbracket$, where $k$ satisfies  $  2 \le k \le m $ since  ${d_n \over D_n}  \underset{n \to + \infty}{\longrightarrow} 0$. Now by using (4.12)-(4.13) and (4.17)-(4.18),  we  can apply the induction hypothesis  above to the $k$ blow-ups: $(x_n^1, r_n^1), ..., (x_n^k, r_n^k)$,  where  $S_n$ is replaced by  $\widetilde{S}_n= 8C_0d_n$, and where $C_0$ is the constant in (4.17).   We obtain  
  $$  \lim_{n \to + \infty} \int_{\bigcup_{j=1}^{k}B_{8C_0d_n}(x_n^j)\setminus  \bigcup_{j=1}^{k}B_{R_n^j}(x_n^j) }e^{4u_n} dV_0 = 0 .  \eqno (4.19) $$

  \bigskip

    \noindent On the other hand, for each fixed $i \in \llbracket 1, k\rrbracket$, if we apply again the induction hypothesis to the  $(m+ 2 -k)$ blow-ups : $x_n^i,  x_n^{k+1}, x_n^{k+2}, ..., x_n^{m+1}$  (we recall here that $2 \le k \le m$) where $R_n^i$ is replaced by $\widetilde{R}_n^i= 8C_0d_n$, then we have for any $ i \in   \llbracket 1, k\rrbracket$, 
    
      $$  \lim_{n \to + \infty} \int_{\left( \bigcup_{j=k +1}^{m+1}B_{S_n}(x_n^j)\bigcup B_{{S}_n}(x_n^i)\right)\setminus \left(  \bigcup_{j=k+1}^{m+1}B_{R_n^j}(x_n^j)\bigcup B_{8C_0d_n}(x_n^i)\right)}e^{4u_n} dV_0 = 0   $$
      which gives 
            $$  \lim_{n \to + \infty} \int_{ \bigcup_{j=1}^{m+1}B_{S_n}(x_n^j)\setminus \left(  \bigcup_{j=k+1}^{m+1}B_{R_n^j}(x_n^j)\bigcup \bigcup_{i=1}^k B_{8C_0d_n}(x_n^i)\right)}e^{4u_n} dV_0 = 0  .  \eqno (4.20) $$
            
            \medskip
            
            \noindent Now it is clear that (4.19) and (4.20) imply 
    
      $$  \lim_{n \to + \infty} \int_{\bigcup_{j=1}^{m+1}B_{S_n}(x_n^j)\setminus\bigcup_{j=1}^{m+1}B_{R_n^j}(x_n^j) }e^{4u_n} dV_0 = 0 .  $$
 This achieves the proof of Proposition 4.3. 
    \end{proof}
  
  \bigskip
  
   The following proposition is the principal tool in the proof of Theorem 1.1 
   
   \medskip

 \begin{prop} Let $(u_n, f_n)$ as in Theorem 1.1.  If the first alternative in Theorem 1.1 does not hold, then there is exist a finite number of blow-ups $(x_n^1, r_n^1)_n , ...., (x_n^k, r_n^k)_n$ with $1\le k \le  {k_0 \over 16\pi^2 }$,  and $k$ sequences  $(R_n^1)_n, ..., (R_n^k)$  of  positive numbers  such that  
 $$ \lim_{n\to + \infty}  R_n^i= 0   \  \  , \ \  \lim_{n\to + \infty}{r_n^i \over  R_n^i} = 0  \  \  \forall \ i \in \llbracket 1, k \rrbracket ,  \eqno (4.21) $$
    $$ \lim_{n\to + \infty} {R_n^i \over  \displaystyle \underset{j\not= i}{\inf_{ 1\le j \le k}} |x_n^i -x_n^j |}   =  0 , \   \   \forall \ i  \in  \llbracket 1, k \rrbracket \ \    \hbox{if } \ k \ge 2,     \eqno (4.22) $$
    and 
 $$ \forall \ B_r(y) \subset  M\setminus \bigcup_{j=1}^kB_{{1\over 2}R_n^j}(x_n^j), \   \int_{B_r(y)} e^{4u_n} dV_0 \ge {\pi^2 \over k_0}  \  \Longrightarrow \   r \ge \alpha \hspace{0,5mm} 
 d_n(y) , \eqno (4.23) $$
 where $\alpha $ is a positive constant   independent of $n$, and where   $ \displaystyle d_n(y) = \inf_{1\le j\le k} |y-x_n^j| $. 
Moreover  we have, for all $i \in \llbracket 1, k \rrbracket$, 
   $$\lim_{n\to + \infty} \int_{B_{R_n^i}(x_n^i)} e^{4u_n} dV_0 =  {16 \pi^2 \over k_0} .  \eqno (4.24)$$
   
 \end{prop}

 \bigskip
 
 \begin{proof}  Before proving Proposition 4.4 let us introduce some notations. If $(x_n^1, r_n^1)_n , ...., (x_n^l , r_n^l)_n$ are $l$  blow-ups for $(u_n)_n$, we say that they satisfy the property $( {\mathcal P} ) $ if 
$$  \displaystyle l = 1 \  \  \hbox{or} \  \    \lim_{n\to + \infty}{r_n^i \over  \displaystyle \underset{j\not= i}{\inf_{ 1\le j \le l}} |x_n^i -x_n^j |}   =   0 , \   \   \forall \ i  \in  \llbracket 1, l \rrbracket \ \    \hbox{if } \ l \ge 2.$$

\medskip

 Now let us prove Proposition 4.4. As noted in the begining of this section,we may suppose that  $k_0 >0$. If the first alternative in Theorem 1.1 does not hold, then by using Proposition 2.3,    there exists a point $x \in M$ such that  for any $r > 0$ we have 
 $$\liminf_{n\to +\infty} \int_{B_r(x)} e^{4u_n} dV_0 \ge {8 \pi^2 \over k_0} + o_r(1).  $$
 where $o_r(1) \to 0$ as $r \to 0$. 
 It  follows that there exit   $x_n \in M$ and $r_n >0$ such that 
 $${ \pi^2 \over k_0} = \int_{B_{r_n}(x_n)} e^{4u_n} dV_0  = \sup_{ x \in M}\int_{B_{r_n}(x)} e^{4u_n} dV_0 \  \  \hbox{and} \  \   \lim_{n\to + \infty}r_n = 0. $$
 
 \noindent Then  setting $\widehat{r}_n = \sqrt{ r_n}$ we have that for any $y \in  B_{\widehat{r}_n}(x_n)$,  \ $\int_{B_{r_n}(y)} e^{4u_n} dV_0 < {\pi^2 \over k_0}$,  and applying Proposition 2.5, we see that $(x_n, r_n)$ is a blow-up for $(u_n)_n$. 
 It follows that the set $A$ defined by 
 
 $$A:= \left\{ l \in \mathbb{N}  \  :  \ \hbox{there exist} \ l \ \hbox{blow-ups} \  (x_n^1, r_n^1), ..., (x_n^l, r_n^l)  \ \hbox{satisfying the property $( {\mathcal P} )$ }   \right\} $$
 is not empty.  First we shall prove that  if $l \in A$, then $l \le {16\pi^2 \over k_0}\ .$ Indeed, let $ l \in A$.  Then  there exit $l$ blow-ups  $(x_n^1, r_n^1), ..., (x_n^l, r_n^l) $ satisfying the property $( {\mathcal P} )$ above.  More precisely, we have 
 $$\lim_{n\to + \infty} {r_n^i \over d_n^i } = 0 \  \  \forall \ i \in \llbracket 1, l \rrbracket ,   \eqno (4.25)$$
 where $\displaystyle d_n^i = \underset{j\not= i}{\inf_{ 1\le j \le l}} |x_n^i - x_n^j| $  if $ l \ge 2$,  and  $d_n^1 = 1$ if $l= 1$.

Now we apply  Proposition 2.3 by  setting   $\beta_n = {d_n^i \over 4 r_n^i}$. Then  there exists  $(b_n^i)_n$ satisfying   $b_n^i \le \beta_n$ and  $b_n^i \underset{n \to + \infty}{\longrightarrow} + \infty $,  such that 
 $$\lim_{n\to + \infty} \int_{B_{b_n^ir_n^i}(x_n^i)} e^{4u_n} dV_0 =  {16 \pi^2 \over k_0}.  \eqno(4.26)$$
 Since $b_n^ir_n^i \le {1\over 4} d_n^i$, then  the balls $B_{b_n^ir_n^i}(x_n^i), \  i= 1, ..., l, $ are pairwise disjoint.  This  implies by using (4.26) that 
 $${16 \pi^2 l \over k_0} = \lim_{n\to + \infty} \int_{\bigcup_{i=1}^lB_{b_n^ir_n^i}(x_n^i)} e^{4u_n} dV_0 \le \int_M e^{4u_n} dV_0 = 1,$$
 which implies $l \le {k_0 \over 16 \pi^2}$.  Hence  the set $A$ defined above is bounded, so let  $k := \max A$. Thus, there exit $k$ blow-ups $(x_n^1, r_n^1)_n , ...., (x_n^k, r_n^k)_n$ satisfying the property $({\mathcal P})$ defined above. That is, 
$$\lim_{n\to + \infty}{r_n^i \over d_n^i } = 0 \   \   \forall \ i \in   \llbracket 1, k \rrbracket ,  \eqno (4.27) $$

\medskip

\noindent where $\displaystyle d_n^i = \underset{j\not= i}{\inf_{ 1\le j \le l}} |x_n^i - x_n^j| $  if $ k \ge 2$,  and  $d_n^1 = 1$ if $k= 1$. 
Now, by setting $\beta_n = {1\over 2}\sqrt{{d_n^i \over  r_n^i}}$  and applying Proposition 2.3, then  there exists  $(b_n^i)_n$ satisfying   $b_n^i \le {1\over 2}\sqrt{{d_n^i \over  r_n^i}} $ and  $b_n^i \underset{n \to + \infty}{\longrightarrow} + \infty $,  such that 
 $$    \lim_{n\to + \infty} \int_{B_{b_n^ir_n^i}(x_n^i)} e^{4u_n} dV_0 =  {16 \pi^2 \over k_0} ,  \  \   \forall i = 1, ..., k.    \eqno(4.28)$$
 If we set $R_n^i =  b_n^ir_n^i$, then  it is clear that (4.21) is satisfied, and (4.24) follows from (4.28). If we apply again Proposition 2.3 by choosing $\beta_n = {1\over 4}{R_n^i \over r_n^i}$, and using (4.28) we arrive at 
  $$\lim_{n\to +\infty} \int_{B_{R_n^i}(x_n^i)\setminus B_{{1\over 4}R_n^i}(x_n^i)} e^{4u_n} dV_0 = 0,  \  \   \forall i = 1, ..., k.   \eqno (4.29) $$
  Hence 
 $$\lim_{n\to + \infty} \int_{B_{{1\over 4}R_n^i}(x_n^i)} e^{4u_n} dV_0 =  \lim_{n\to + \infty} \int_{B_{R_n^i}(x_n^i)} e^{4u_n} dV_0  = {16 \pi^2 \over k_0} ,  \  \   \forall i = 1, ..., k.    \eqno(4.30)$$
 
    \medskip

   Now, Since $R_n^i  = b_n^ir_n^i \le \sqrt{r_n^i d_n^i}$, then we have 
  $$\lim_{n\to + \infty} {R_n^i \over d_n^i } = 0 \  \  \forall \ i \in \llbracket 1, k \rrbracket ,   \eqno (4.31)$$
  which proves (4.22). 
  
  \medskip

  It remains then  to prove (4.23). Suppose by contradiction that (4.23) is false, then there are balls $ \displaystyle B_{\rho_n}(z_n) \subset M\setminus \bigcup_{j=1}^kB_{{1\over 2}R_n^j}(x_n^j)$ such that 
 $$\int_{B_{\rho_n}(z_n)} e^{4u_n} dV_0 \ge {\pi^2 \over k_0}  \   \  \hbox{and} \  \   \lim_{n\to +\infty} {\rho_n \over d_n(z_n)} = 0,$$
 \medskip
 
\noindent  where we recall that  \ $\displaystyle d_n(z) = \inf_{1 \le j\le k}|x_n^j - z|$.  Then there exist $r_n \le \rho_n$ and a ball $B_{r_n}(a_n) \subset M\setminus \bigcup_{j=1}^kB_{{1\over 2}R_n^j}(x_n^j)$ such that 
$$ { \pi^2 \over k_0} = \int_{B_{r_n}(a_n)} e^{4u_n} dV_0  = \sup_{B_{r_n}(y) \subset M\setminus \bigcup_{j=1}^kB_{{1\over 2}R_n^j}(x_n^j)}\int_{B_{r_n}(y)} e^{4u_n} dV_0. \eqno (4.32)$$
Let us show  that 
$$\lim_{n \to + \infty}{r_n \over d_n(a_n)} = 0 . \eqno (4.33)$$

If (4.33) was false, then by passing to a subsequence if necessary,  we would have for some constant $C$ independent of $n$,
$$ r_n \ge C d_n(a_n)   \eqno (4.34)$$
and  without loss of generality we may suppose that $d_n(a_n) = |x_n^1 - a_n|$.  Set $d_n := d_n(a_n) $  and  define :
 $$J :=  \{ \  j \in \llbracket 1, k \rrbracket  \ : \  |x_n^j-x_n^1| \le Cd_n   \   \forall n \  \}  , $$
 
 \medskip

 \noindent where $C$ is a positive constant independent of $n$. 
  Observing that $ 1 \in J$, so by relabeling the blow-ups, we may suppose  that  $J= \llbracket 1, m\rrbracket$, with $  1 \le m \le k $, and by 
  passing to a subsequence if necessary,  we have 
     $$  \forall j \in 
  \llbracket 1, m \rrbracket , \   |x_n^j -x_n^1 |\le C_0 d_n  \   \   \forall n  ,  \eqno (4.35)$$
  and
   $$  \forall j \in 
  \llbracket m+1, k \rrbracket , \  \lim_{n\to + \infty}{d_n  \over  |x_n^j -x_n^1 | } = 0 ,  \eqno (4.36) $$
  
    \medskip
    
\noindent where  $C_0$ is a positive constant independent of   $n$ that we assume,  without loss of generality,  satisfying $ C_0 \ge 1$.

  \medskip

Now by using  (4.32), (4.34) and (4.35) one  can easily check that
$$ \forall \ B_r(y) \subset  \bigcup_{j=1}^mB_{8C_0d_n}(x_n^j)\setminus \bigcup_{j=1}^mB_{{1\over 2}R_n^j}(x_n^j), \   \int_{B_r(y)} e^{4u_n} dV_0 \ge { \pi^2 \over k_0}  \  \Longrightarrow \   r \ge \alpha \hspace{0,5mm} 
 \inf_{1 \le j \le m}|y-x_n^j|, \eqno (4.37) $$
 where $\alpha $ is a positive constant independent of   $n$. Thus by applying Proposition 4.3, where we take $S_n= 4C_0d_n$ and using (4.37), we get 
 
     $$ \lim_{n \to + \infty} \int_{B_{4C_0d_n}(x_n^1)\setminus\bigcup_{j=1}^{m}B_{R_n^j}(x_n^j)}e^{4u_n} dV_0 = 0 ,$$
     
   \noindent  which contradicts (4.32) since  $B_{r_n}(a_n) \subset  B_{4C_0d_n}(x_n^1)\setminus\bigcup_{j=1}^{m}B_{R_n^j}(x_n^j)$. So this proves (4.33). 
 
 \medskip

Since  ${1\over 2} R_n^i \le | x_n^i - a_n |$, then  it follows from  (4.29), (4.32) and (4.33) that for $n$ large enough we have 
 $$R_n^i \le {4 \over 3}| x_n^i - a_n |  \  \  \forall i = 1, ..., k.  \eqno (4.38) $$

\medskip
 \noindent Indeed, if (4.38)  were  not satisfied, then by passing to a subsequence one could check   by using (4.33) that $B_{r_n}(a_n) \subset B_{R^i_n}(x^i_n) \setminus B_{{1\over 4}R^i_n}(x^i_n)$, so by (4.29) we would have that $\displaystyle \lim_{n \to +\infty}\int_{B_{r_n}(a_n) }e^{4u_n} dV_0 = 0 $ contradicting thus (4.32). 
 
 \medskip

 Now,  by  using Proposition 2.5, where we take  $x_n = a_n$, $\widehat{r}_n = {1 \over 4} d_n$,  and using (4.33) and (4.38), it is not difficult to see that     $(a_n, r_n)$ is a blow-up for $(u_n)_n$,  and by using (4.27) and  (4.33) we see  that  the $(k+1)$ blow-ups $(x_n^1, r_n^1), ..., (x_n^k, r_n^k), (a_n, r_n)$ satisfy the property $({\mathcal P})$. This contradicts the fact that $k = \max A$. The proof of Proposition 4.4 is then complete.

 \end{proof}
 
 \bigskip
 
 Now we are in position to prove Theorem 1.1. 
 
 \medskip

 \begin{proof}[Proof of Theorem 1.1]
 
 \medskip
 
 Let $(u_n, f_n)$ as in Theorem 1.1. If the first alternative in Theorem 1.1 does not hold, then by Proposition 4.4 there are $k$ blow-ups $ (x_n^1, r_n^1)_n , ...., (x_n^k, r_n^k)_n$ with $1\le k \le  {k_0 \over 16\pi^2 }$,  and $k$ sequences  $(R_n^1)_n, ..., (R_n^k)_n$  of  positive numbers  satisfying (4.21)-(4.24) in Proposition 4.4.  If we apply Proposition 4.3 by taking \  $S_n = 2  \hbox{ diam}(M)$, we obtain 
 $$\lim_{n \to + \infty}\int_{M\setminus \bigcup_{i=1}^kB_{R_n^i}(x_n^i) } e^{4u_n} dV_0 = 0 ,  \eqno (4.39) $$
which implies since the balls $B_{R_n^i}(x_n^i)$ are pairwise disjoint,  
 $$  \lim_{n \to +\infty}\int_M e^{4u_n} dV_0  =  \lim_{n \to + \infty} \sum_{i=1}^k \int_{B_{R_n^i}(x_n^i)}e^{4u_n} dV_0 .  \eqno (4.40) $$
 
 \medskip
 
  Since by (4.24) we have 
 $$\lim_{n\to +\infty} \int_{B_{R_n^i}(x_n^i)}e^{4u_n} dV_0  = {16\pi^2 \over k_0}  \  \  \forall i = 1, ..., k, \eqno (4.41) $$
 and since $\displaystyle  \int_M e^{4u_n} dV_0 = 1 $, then we get from (4.40) and (4.41) that
  $$k = {k_0 \over 16 \pi^2} \ . \eqno (4.42) $$
 
 \medskip

 On the other hand,  since $M$ is compact, then  by passing to  a subsequence,  there exist   $m$ distincts points $a_1, ..., a_m \in M$  with $m \le k$ such that for any $i = 1, ..., k$, the sequence $(x_n^i)_n$ converges to a limit in $\{ a_1, ..., a_m \}$.  For any $i= 1, ..., m$, if we set 
 $$ l_i = \#\{  \ j \in  \llbracket 1, k \rrbracket \ : \ \lim_{n\to \infty} x_n^j = a_i \ \} , \eqno (4.43) $$
 then we have 
 $$l_1 + \cdots + l_m = k =  {k_0 \over 16 \pi^2} \  ,   \eqno (4.44) $$
 where we have used (4.42). 
 
 \medskip
 
 Let now $\varphi \in C^{0}(M)$. Then we have by (4.39) 
$$\lim_{n \to +\infty}  \int_M\varphi e^{4u_n}dV_0 = \lim_{n \to + \infty} \sum_{i=1}^k \int_{B_{R_n^i}(x_n^i)}\varphi e^{4u_n} dV_0 .  \eqno (4.45)$$
But we have by the mean-value Theorem 
$$ \int_{B_{R_n^i}(x_n^i)}\varphi e^{4u_n} dV_0 = \varphi(y_n^i) \int_{B_{R_n^i}(x_n^i)}e^{4u_n} dV_0 \eqno (4.46) $$
for some $y_n^i \in B_{R_n^i}(x_n^i)$. Since $R_n^i \to 0$ as $n \to \infty$, then we have 
$$\displaystyle \lim_{n \to + \infty}y_n^i = \lim_{n \to + \infty}x_n^i  \in \{ \ a_1, ..., a_m \ \} . \eqno (4.47) $$

\medskip

\noindent It follows from (4.45) by using  (4.43), (4.44), (4.46) and (4.47) that 
$$\lim_{n \to +\infty}  \int_M\varphi e^{4u_n}dV_0  ={16 \pi^2\over k_0}  \sum_{i=1}^m l_i  \varphi(a_i) .$$
This achieves the proof of Theorem 1.1. 
 
 \end{proof}
 
 \bigskip

\section{The flow} 

\medskip

In this section we prove our results concerning the $Q$-curvature flow.   Through this section we assume that the total $Q$-curvature $k_0$ satisfies  $k_0 > 0$ since  $k_0\le 0 $ is included in the case $k_0 \le 16\pi^2$ which has been already proved by S. Brendle \cite{sB1}. 

\bigskip

 \newtheorem{lem}{Lemma}[section]

\begin{lem} Let  $u \in C^{\infty}(M\times [0, T))$   be the solution of problem  $(1.15)$  defined on a maximal interval $[0, T)$, and set  
$$A_t := \bigl\{ \ x \in M \ : \  u(t,x) \ge \alpha_0 \ \bigr\} ,  \   \   t \in [0, T ), $$
where $ \displaystyle \alpha_0 ={1\over 4}\log\left( {1\over 2|M|}\int_M e^{4u_0} dV_0\right)$, and where $|M|$ is the volume of $(M, g_0)$. For any $L_0 > 0$,  there exists a positive constant $C_0$ depending only  on $L_0$ and $M$ such that,  for any $T_0 \in [0, T)$, if 

 $$\| u_0\|_{H^2(M)} \le L_0  \hspace{3mm} \hbox{and} \hspace{2mm}  \inf_{t\in [0, T_0]} E(u(t)) \ge - L_0,   \eqno (5.1) $$
 
  \noindent then   $A_t$ has volume $|A_t|$ $($with respect to $g_0$$)$  satisfying 
  
$$|A_t| \ge  \exp\left(- C_0e^{ 2 k_0 T_0}\right) \  \   \hbox{for all}  \    t \in [0, T_0].  \eqno (5.2) $$

\end{lem} 

\medskip

\begin{proof} 
 Through the proof of Lemma 5.1, $C$ will denote a positive constant depending only  on $L_0$ and $M$, whose value may change from line to line.

\medskip

 Since by (1.16) the volume of the conformal metric $e^{2u(t)}g_0$ remains constant, we may assume without loss of generality that for all $t \in [0, T)$ 
$$\int_M e^{4u(t)} dV_0 = 1.  \eqno (5.3) $$

\noindent Thus the first equation in (1.15) becomes 
 $$e^{4u}\partial_t u  = -{1 \over 2}\left( P_0u + Q_0 \right) + {1 \over 2}k_0e^{4u} .  \eqno (5.4) $$
  \noindent Multiplying equation  (5.4)  by $u(t)$ and integrating on $M$ with respect to $dV_0$,  and using (5.3),  one gets 
 $$ {d \over dt} \int_M ue^{4u} \ dV_0 =   -2 \int_M P_0 u\cdot u  \ dV_0  - 2  \int_M Q_0 u \ dV_0 + 2k_0 \int_M ue^{4u} \ dV_0 . \eqno (5.5) $$
 
 \medskip
 
 \noindent Let $L_0 \in \Ree$ and $T_0 \in [0, T)$ such  that (5.1) is satisfied. Then we have  for any $t\in [0, T_0]$
 $${1\over 2}\int_M P_0 u\cdot u  \ dV_0 +  \int_M Q_0 u \ dV_0  = E(u(t))  \ge   - L_0 .   \eqno (5.6)$$

 \medskip
 
 \noindent   It follows from (5.5) and (5.6) that 
 $${d \over dt} \int_M ue^{4u} \ dV_0  \le - \int_M P_0 u\cdot u  \ dV_0 + 2k_0  \int_M ue^{4u} \ dV_0 + 2 L_0 $$
 
\noindent which implies since $P_0$ is positive

  $${d \over dt} \int_M ue^{4u} \ dV_0  \le 2 k_0  \int_M ue^{4u} \ dV_0 + 2  L_0.  \eqno (5.7)$$
  
  \medskip
  
  \noindent By setting $\displaystyle Y(t) = \int_M ue^{4u} \ dV_0$, it follows from (5.7) that  for all $t \in [0, T_0]$,  
  $$Y(t)  \le  \left( Y(0)  + { L_0 \over k_0}\right) e^{2k_0t}  \le Ce^{2 k_0T_0}, \eqno (5.8)$$
  where the constant $C$ depends only  $L_0$ and $M$  since $Y(0)$ depends only on the $H^2$-norm of $u_0$ by Adams inequality (see section 2). 
  
  \medskip

   Since $ue^{4u} \ge -{e^{-1}\over 4 }$, then  we get from (5.8),  for any $A \subset M$, 
  $$\int_A ue^{4u } \ dV_0 \le  Ce^{2k_0T_0} . \eqno (5.9) $$
  
  \medskip
  
  For $z>0$, let $\varphi(z) = z\log z$. Then $\varphi$ is convex on $(0, + \infty)$, and it satisfies for each $\lambda > 1$ and $z> 0$, 
  $$z = {\varphi(\lambda z) \over \varphi(\lambda)} - {\varphi(z) \over \log\lambda} , $$
  
  \noindent which implies, since $\varphi(z) \ge - e^{-1} $  for any $z > 0$,
    $$z \le  {\varphi(\lambda z) \over \varphi(\lambda)} + {e^{-1}\over  \log\lambda} . \eqno (5.10) $$
    
    \medskip
    
  For $t \in [0, T_0]$,  let $A_t \subset M$ defined by 
   $$ A_t = \{ \ x \in M \ : \  u(x, t) \ge \alpha_0 \ \} $$
   where 
   $$\alpha_0 = {1\over 4}\log\left( {1\over 2|M|}\int_M e^{4u_0} dV_0\right) = {1\over 4}\log\left( {1\over 2|M|}\right) \  .$$
   
   \medskip
   
    \noindent Since $\varphi$ is convex, then it follows from Jensen inequality 
    
    $$\varphi\left({1 \over |A_t|}\int_{A_t} e^{4u} dV_0 \right) \le {1 \over |A_t|}\int_{A_t} \varphi\left(e^{4u}\right) dV_0  \eqno (5.11) $$
   
   \noindent But by (5.9) we have 
   $$ {1 \over |A_t|}\int_{A_t} \varphi\left(e^{4u}\right) dV_0  \le  {Ce^{2k_0T_0} \over |A_t|}, $$
   hence it follows from (5.11) that 
   $$\varphi\left({1 \over |A_t|}\int_{A_t} e^{4u} dV_0 \right) \le {Ce^{2k_0T_0}\over |A_t|} . \eqno (5.12)$$
   
   \medskip
   
   \noindent Now, if $|A_t| \ge 1$, then the estimate (5.2) is trivially satisfied by taking $C_0$ any positive constant, and  Lemma 5.1  is proved in this case. Thus we may suppose that $ |A_t| < 1$.  Then by using (5.10) with $\displaystyle \lambda = {1 \over |A_t|}$ and $ \displaystyle z = \int_{A_t} e^{4u} dV_0$, we have 
   $$\int_{A_t} e^{4u} dV_0 \le {|A_t| \over \log{1 \over |A_t|}}  \ \varphi\left({1 \over |A_t|}\int_{A_t} e^{4u} dV_0 \right) + {e^{-1} \over  \log{1 \over|A_t|}} \ ,$$ 
   which gives by using (5.12)
   $$\int_{A_t} e^{4u} dV_0  \le  \left( Ce^{2 k_0T_0} + e^{-1} \right) {1 \over \log{1 \over |A_t|}} \le  { Ce^{2 k_0T_0}  \over \log{1 \over |A_t|}}       \  . \eqno (5.13) $$
   
   \medskip
   
   \noindent On the other hand, we have 
   $$1 = \int_M e^{4u} dv_0 = \int_{A_t} e^{4u} dV_0 + \int_{M\setminus A_t} e^{4u} dV_0  \eqno (5.14)  $$
   and since $e^{4u}  <e^{4 \alpha_0} =   {1 \over 2 |M|}$ on $M\setminus A_t$, \  then (5.14) implies
   $$ {1 \over 2} \le \int_{A_t} e^{4u} dV_0  $$
   which together with (5.13) give 
   $$ \log{1 \over |A_t|} \le  Ce^{2 k_0T_0}. $$
   
   \medskip
   
  \noindent  This achieves the proof  Lemma 5.1. 
 \end{proof}

 \bigskip

 Lemma 5.1 allows us to prove the following estimates  on the solution : 
 
 \medskip
 
 \begin{prop} Let  $u \in C^{\infty}(M\times [0, T)$   be the solution of problem $(1.15)$ defined on a maximal interval $[0, T)$. For any $L_0 > 0$,  there exists a positive constant $C_0$ depending on $L_0$ and $M$ such that,  for any $T_0 \in [0, T)$, if 
 
 $$\| u_0\|_{H^2(M)} \le L_0  \hspace{3mm} \hbox{and} \hspace{2mm}  \inf_{t\in [0, T_0]} E(u(t)) \ge - L_0,  $$
then we have 
$$\sup_{t \in [0, T_0]} \|u(t)\|_{H^2(M)} \le \exp\left(C_0e^{2k_0 T_0}\right) . \eqno (5.15) $$
Moreover, for any $ k \in \mathbb{N}$, there exist a positive constant $C_k$ depending on $k, L_0 , T_0 $ and $M$ such that 
$$\sup_{t \in [0, T_0]} \|u(t)\|_{H^k(M)} \le C_k . \eqno (5.16) $$
\end{prop} 
 
 \medskip
 
 \begin{proof}  Through the proof of Proposition 5.1, $C$ will denote a positive constant depending only  on $L_0$ and $M$, whose value may change from line to line. For any measurable set $A \subset M$, we shall denote its volume with respect to the metric $g_0$  by $|A|$.

Since by (1.16) the volume of the conformal metric $e^{2u(t)}g_0$ remains constant, we may assume without loss of generality that 
$$\int_M e^{4u(t)} dV_0 = 1.  \eqno (5.17) $$
 
 \noindent This implies by using the elementary inequality $z \le e^z$, that 
  $$ \int_A u(t) dV_0 \le {1 \over 4}  \int_A e^{4u(t)} dV_0 \le {1\over 4}  \eqno (5.18)$$
  for any $A \subset M$.   Let $T_0 \in [0, T)$ and $L_0 >0 $  such that 
  $$    \inf_{t\in [0, T_0]} E(u(t)) \ge -  L_0.  $$
  If we let $A= A_t$, where $A_t$ is as in Lemma 5.1, then we have by using (5.18) and the definition of the set $A_t$,  for any $t \in [0, T_0], $ 
   $$ \left|\int_M u(t) dV_0\right|   \le  \left| \int_{A_t }u(t) dV_0\right| +   \left|\int_{M\setminus A_t} u(t) dV_0 \right|  \eqno  $$
  $$\le C  +   \left|\int_{M\setminus A_t} u(t) dV_0 \right| . \eqno (5.19)$$ 
  
    \noindent But by the Cauchy-Schwarz inequality we have 
  $$\left|\int_{M\setminus A_t} u(t) dV_0 \right|  \le |M\setminus A_t|^{1\over 2} \|u\|_{L^2(M)}, $$
  and by replacing this inequality  in (5.19),  we get for any $\varepsilon > 0$,  
  $$   \left(\int_M u(t) dV_0\right)^2   \le   \left(1+ \varepsilon \right)  |M\setminus A_t| \|u(t)\|_{L^2(M)}^2 + C \varepsilon^{-1}  + C . \eqno (5.20) $$

  \medskip
  
  Now, from  Poincar\'e's  inequality we have 
  $$\|u(t)\|_{L^2(M)}^2 \le  {1 \over \lambda_1} \int_M P_0 u(t)\cdot u(t)  \ dV_0 + |M| \  |\overline{u}(t)|^2,  \eqno (5.21)$$

  \noindent where $\lambda_1$ is the first positive eigenvalue of $P_0$,  and   \ $\displaystyle \overline{u}(t) = {1 \over |M|} \int_{M}u(t) dV_0$
  is the average of $u(t)$.  Thus it follows from (5.20) and (5.21) that 
  $$\left( 1 - {(1+ \varepsilon ) |M\setminus A_t| \over |M|} \right) \|u(t)\|_{L^2(M)}^2  \le {1 \over \lambda_1} \int_M P_0 u\cdot u  \ dV_0  + {C\over |M|} \varepsilon^{-1} + {C\over |M|}  $$
  that is 
  $$  \left(| A_t|  - \varepsilon  |M\setminus A_t| \right) \|u(t)\|_{L^2(M)}^2  \le {|M| \over \lambda_1} \int_M P_0 u\cdot u  \ dV_0  + C \varepsilon^{-1} + C.   \eqno (5.22)$$
  
  \medskip
  
  \noindent Since by Lemma 5.1 we have $|A_t|\ge \exp\left(-C_0e^{2 k_0T_0}\right)$, then  by choosing 
  $ \varepsilon = {1\over 2|M|} \exp\left(-C_0e^{2 k_0T_0}\right)$ in (5.22)  and observing that $ |M\setminus A_t| \le |M|$, we obtain 
  $$\|u(t)\|_{L^2(M)}^2 \le   C \left( \int_M P_0 u(t)\cdot u(t)  \ dV_0  +  1 \right)  \exp\left(C_0e^{2 k_0T_0}\right).  \eqno (5.23) $$

\medskip

\noindent Since the functional $E$ is decreasing along the flow by (1.17), then 
$${1\over 2} \int_M P_0 u(t)\cdot u(t)  \ dV_0   + \int_M Q_0 u(t) \ dV_0 = E(u(t)) \le E(u_0), $$
hence
$$  \int_M P_0 u(t)\cdot u(t)  \ dV_0    \le C \|u(t)\|_{L^2(M)} + C. \eqno (5.24) $$

\medskip

\noindent It follows from (5.23) and (5.24) that 
$$\|u\|_{L^2(M)}^2 \le   \left( \|u(t)\|_{L^2(M)} + 1 \right)  \exp\left(Ce^{2 k_0T_0}\right),$$
which implies that 
$$  \|u(t)\|_{L^2(M)}  \le   \exp\left(Ce^{2 k_0T_0}\right) .   \eqno (5.25) $$
 
Combining (5.24) and (5.25) we get (5.15).  The higher order estimate (5.16) follows as in S. Brendle \cite{sB1}.

 \end{proof}

 \bigskip

\begin{proof}[Proof of Theorem 1.2]

{\bf Step 1 \ Global existence of the flow}.   \  Let   $u \in C^{\infty}(M\times [0, T))$   be the solution of problem (1.15) defined on a maximal interval $[0, T)$, satisfying (1.18), that is 
$$L :=  \inf_{t\in [0, T)} E(u(t)) > - \infty .  \eqno (5.26)  $$ 

Suppose by contradiction that $T< + \infty$,   then it follows from Proposition 5.1 by taking  $L_0 =  \|u_0\|_{H^2(M)} + |L|$ \  that 

$$ \sup_{t \in [0, T)}\|u(t)\|_{H^2(M)} <  \exp\left(C_0e^{2 k_0T}\right) ,  $$
and for  for any $k \ge 2$ :
 $$ \sup_{t \in [0, T)}\|u(t)\|_{H^k(M)} < + \infty  .  \eqno (5.27) $$
It is clear that (5.27) implies that   the solution $u(t)$ would be extended beyond $T$ giving thus a contradiction. This proves Step 1.

  \bigskip

\noindent {\bf Step 2 \  Convergence of the flow}.  According to the first step, the solution $u$ is defined on $[0, + \infty)$, and (5.26) becomes 
$$L:= \inf_{t\in [0, + \infty)} E(u(t)) > - \infty .   \eqno(5.28)$$
 Since by (1.16) the volume of the conformal metric $e^{2u(t)}g_0$ remains constant, we may assume without loss of generality that 
$$\int_M e^{4u(t)} dV_0 = 1.  \eqno (5.29) $$
By using (1.17), we get  for any $T> 0$,
$$\int_{0}^{T} \int_Me^{4u(t)} |\partial_tu(t)|^2 dV_0 dt  = E(u_0) - E(u(T)) \ \le  E(u_0)  - L  ,$$
which implies 
$$\int_{0}^{+\infty} \int_M e^{4u(t)} |\partial_tu(t)|^2  dV_0dt   \le  E(u_0)  - L. \eqno (5.30)$$

By using the mean value theorem, we obtain from (5.30), that for any $n \in \Nee$, there exits $t_n \in [n, n+ 1]$ such that 
$$\lim_{n\to +\infty} \int_Me^{4u(t_n)} |\partial_tu(t_n)|^2 dV_0   =  0 . \eqno (5.31) $$
Now if we set 
$$u_n= u(t_n) \   \   \hbox{and} \   \    f_n = 2 e^{4u_n}\partial_tu(t_n) + Q_0$$
then we have  from (1.15) 
$$ P_0 u_n + f_n =  k_0 e^{4u_n} \eqno (5.32) $$
with 
$$ \int_M e^{4u_n} dV_0 = 1$$
and 
$$\|f_n - Q_0\|_{L^1(M)} \le2 \left(\int_M e^{4u_n} dV_0 \right)^{1/2}\left(\int_Me^{4u_n}|\partial_tu(t_n)|^2 dV_0 \right)^{1/2} $$
$$=   2 \left(\int_Me^{4u_n}|\partial_tu(t_n)|^2 dV_0 \right)^{1/2}
 \underset{n \to + \infty}\longrightarrow  0 .$$
Since we are supposing $k_0 \not\in 16\pi^2\Nee^*$, then we can apply Corollary 1.1  to get, for any $p\ge 1$
 $$\int_Me^{p|u_n|} dV_0 \le C_p, \eqno (5.33)$$
 which implies by using (5.31) that for any $q \in [1, 2 )$ 
 $$\|f_n\|_{L^q(M)} \le C_q \eqno (5. 34)$$
 for some constant $C_q$ depending on $q$.  Thus it follows from the elliptic regularity theory applied to equation (5.32) by using (5.33) and (5.34) that $(u_n)_n$ is bounded in $W^{4,q}(M)$ for any $ q \in [1, 2)$.  But by  Sobolev embedding theorem we have  $W^{4,q}(M) \subset C^{\alpha}(M)$ for any $\alpha \in (0, 1)$, and by applying the elliptic regularity theory again to equation (5.32), we obtain that $(u_n)_n$ is bounded in $H^4(M)$. In particular we have that $(u_n)_n$ is bounded  in $H^2(M)$, that is 
 $$\|u_n\|_{H^2(M)} \le C ,\eqno (5.35)$$
 where $C$ is a positive constant depending only on $L, u_0$ and $M$.  Now, let  us define $v_n(t) := u(t+t_n)$. Then $v_n$ is a solution of problem (1.15) where $u_0$ is replaced by $u_n$, that is 
 $$ \displaystyle 
\begin{cases} \partial_t v_n  = -{1 \over 2} e^{-4v_n} \left( P_0 v_n + Q_0 \right) + {k_0 \over 2} \cr  \cr v_n(0)= u_n.   \end{cases}
\eqno (5.36)
$$
 We want to apply Proposition 5.1 to $v_n$.  We  have 
 $$ \inf_{t\in [0, 1]}E(v_n(t)) =  \inf_{t\in [0, 1]}E(u(t + t_n))  =  \inf_{t\in [t_n, t_n+1]}E(u(t)) \ge L , $$
 where $L$ is given by (5.28). Then by Proposition 5.1, where we choose  $T_0 = 1$ and $L_0 = |L| + C$  with $C$ as in (5.35),  there exist a positive constant $C_0$ depending on $L, u_0$ and $M$, such that 
 $$\sup_{t \in [0, 1]} \|v_n(t)\|_{H^2(M)} \le \exp\left(C_0e^{2 k_0}\right), $$
that is 
  $$\sup_{t \in [t_n, t_n + 1]} \|u(t)\|_{H^2(M)} =  \sup_{t \in [0, 1]} \|v_n(t)\|_{H^2(M)}   \le \exp\left(C_0e^{2 k_0 }\right), $$
  
  \medskip
  
\noindent   and since $n \le t_n \le n+1$ for all $n \in \Nee$, then we have 
  $$\sup_{t \in [0, + \infty)} \|u(t)\|_{H^2(M)}  \le  \exp\left(C_0e^{2 k_0}\right)  . \eqno (5.37)$$
  Following the argument of S. Brendle\cite{sB1}, one gets from (5.37) that, 
    $$\sup_{t \in [0, + \infty)} \|u(t)\|_{H^k(M)}  \le C_k  $$
    for any $k \ge 2$,  and the convergence of the flow follows as in S. Brendle \cite{sB1}. 
\end{proof}

\bigskip

\begin{proof}[Proof of Theorem 1.3]

We proceed by contradiction.   For $u \in C^{\infty}(M)$,  let $\Phi(t,u)$ be the solution of (1.15)    such that $\Phi(0, u) = u$,   that is, 
 $$ \displaystyle 
\begin{cases} \partial_t \Phi  = -{1 \over 2} e^{-4\Phi} \left( P_0\Phi + Q_0 \right) + {1 \over 2}\frac{ k_0 }{\int_M  e^{4\Phi} \ dV_0 } 
 \cr  \cr \Phi(0, u)= u  .  \end{cases}
\eqno (5.38)
$$
Let $[0, T_u)$  be the maximal existence interval of $\Phi$ and suppose by contradiction that 

$$\inf_{t\in [0, T_u)}E(\Phi(t, u)) = -\infty \  \  \forall u \in C^{\infty}(M) . \eqno (5.39) $$

\medskip

Let $X:= C^{\infty}(M)$ endowed with its natural  $C^{\infty}$ topology , and let us introduce the sub-level set 

$$X_0 := \{ \ u \in X \ : \  E(u) \le -L  \ \} , \eqno (5.40) $$

\medskip

\noindent where $L > 0$ is large enough.  One fundamental property of $X_0$ is its invariance under the flow $\Phi$, that is,  if $u \in X_0$, then $ \Phi(t, u) \in X_0$ for all $t \in [0, T_u)$, as it  can be immediately checked by using the fact that $E$ is decreasing along the flow  $\Phi$(see formula 1.17)) 

\medskip

Following  Z. Djadli and A. Machioldi \cite{zD}, one can  prove that   $X_0$  is not contractible. Indeed, in \cite{zD} the set $X_0$ consists of $H^2$- functions $u$ satisfying $E(u) \le - L$, but by following the same proof as in \cite{zD}, one can easily see that the same arguments work when considering $C^{\infty}$ functions and the $C^{\infty}$ topology on $X_0$. Then we shall use our  flow  $\Phi$ to construct a deformation retraction  from $X$ onto $X_0$,  which would give a contradiction since $X$ is contractible as a topological vector space. 

\medskip

By using (5.39)  we can define for any $u \in X$

$$ t_u = \min\{ \ t \in [0, T_u) \ : \  E(\Phi(t, u)) \le - L \ \} . \eqno (5.41)$$

\medskip

\noindent Thus we have  by using the continuity of $\Phi$ that 

$$E(\Phi(t_u, u)) = -L .  \eqno (5.42) $$

\medskip

We extend $\Phi$  on $[0, +\infty)$ by considering  $\widehat{\Phi} : [0, +\infty) \times X  \to X$ as follows 
$$ \widehat{\Phi}(t, u) = \begin{cases}   \Phi(t, u) \  \ \    \hbox{if} \ \  t \in [0, t_u] \cr \cr 
 \Phi(t_u, u) \  \  \hbox{if} \ \  t \ge t_u .
\end{cases}$$
By using Proposition 5.1  and the fact that the functional   $E$  is  decreasing along the flow $\Phi$ (see  formula (1.17)), one can  prove that $\widehat{\Phi}$ is continuous on $[0, +\infty)\times X$. 
\medskip

We define now the following homotopy map : $ H : [0, 1]\times X \to X$ by 
$$ H(t, u) = \begin{cases}   \widehat{\Phi}({t\over 1-t} , u) \  \  \hbox{if} \ \  t \in [0,  1) \cr \cr 
 \widehat{\Phi}(t_u, u) \hspace{5mm}  \hbox{if} \ \     t = 1.
\end{cases}$$

\medskip

\noindent Then it is easy to see that we have 
$$ \begin{cases}  H(0, u) = u    \  \   \forall \  u \in X,   \cr 
H(t , u) = u   \ \   \forall \ u \in X_0, \  \forall \ t \in [0, 1]  \cr
H(1,u) \in X_0 \  \  \forall \ u \in X. 
\end{cases} $$

\medskip

\noindent This proves that $X_0$ is a deformation retract of $X$ which is impossible since $X_0$ is non contractible. The proof of Theorem 1.3 is then complete. 

\end{proof}

\bigskip

\begin{proof}[Proof of Theorem 1.4] Let ${\mathcal S}$ be the set of all solutions of the $Q$-curvature equation 

$$P_0u + Q_0 = k_0 e^{4u}   \eqno (5.43)$$
such that 
$$\int_M e^{4u} dV_0 = 1 $$
(we note here that  by  using (5.43), the last condition is automatically satisfied when $k_0\not= 0$).

\bigskip

According to Corollary 1.1, we have  for any $k \in \mathbb{N}$ 

$$\|u \|_{C^k(M)} \le C_k \ \  \forall u \in { \mathcal S } \eqno (5.44)$$

\medskip

\noindent where $C_k$ is a positive constant independent of $u$.  It follows from (5.44) that the functional $E$ satisfies 

$$E(u) \ge  L   \ \  \forall u \in { \mathcal S } \eqno (5.45)$$
for some constant $L \in \mathbb{R}$ independent of $u$. 

\medskip

Let $\lambda < L$, we shall prove that for any $u_0 \in C^{\infty}(M)$ with $E(u_0) \le \lambda$,  the solution $u(t)$ of (1.15) such that $u(0)= u_0$,  satisfies 
$$\lim_{t\to T}E(u(t)) = - \infty , \eqno (5.46)$$
where $[0, T)$ is the maximal existence interval of $u$. Indeed, suppose by contradiction that (5.46) does not hold. Then according to Theorem 1.2, we have $T= + \infty$, and the solution $u(t)$ converges (as $t\to +\infty$)  to a function $u_{\infty} \in C^{\infty}(M)$ satisfying 
$$ P_0 u_{\infty} + Q_0 ={ k_0 \over \int_Me^{4u_{\infty}} dV_0}  e^{4u_{\infty}}.   \eqno (5.47) $$
Moreover, since the functional $E$ is decreasing along the flow,  $u_{\infty}$ satisfies 
$$E(u_{\infty}) \le E(u_0) \le \lambda .   \eqno (5.48)  $$

\medskip

On the other hand, since $E$ is translation invariant, that is, $E(u+c) = E(u)\  \forall c \in  \mathbb{R}$, we may assume by adding an appropriate constant to $u_{\infty}$, that $ \displaystyle \int_Me^{4u_{\infty}} dV_0 = 1$. This implies by using  (5.47) that 
 \ $u_{\infty} \in {\mathcal S}$. Thus it follows from (5.45) that 

$$E(u_{\infty}) \ge  L  $$

\medskip

\noindent which contradicts (5.48) since $\lambda < L$. This achieves the proof of Theorem 1.4. 

\end{proof}

\end{document}